\documentclass[12pt]{article}
\usepackage{amstext}
\usepackage{amsfonts}
\usepackage{amssymb}
\usepackage{amsbsy}
\usepackage{latexsym}
\usepackage{xy}
\usepackage{hhline}
\xyoption{all}
\newcommand{\vtx}[1]{*+[o][F-]{\scriptscriptstyle #1}} 
\newtheorem{theo}{Theorem}

\newtheorem{cor}{Corollary}
\newtheorem{lemma}{Lemma}
\newtheorem{remark}{Remark}

\newcommand{\Ref}[1]{(\ref{#1})}

\newcounter{theexample}
\newcommand{\example}
{\refstepcounter{theexample}%
\bigskip\noindent{\bf Example~\arabic{theexample}.}}

\newenvironment{proof}{{\it Proof. }}
{$\Box$ \bigskip}

\newenvironment{eq}{\begin{equation}}{\end{equation}}

\newcommand{\si}{\sigma}
\newcommand{\al}{\alpha}
\newcommand{\be}{\beta}
\newcommand{\ga}{\gamma}
\newcommand{\la}{\lambda}
\newcommand{\de}{\delta}
\newcommand{\De}{\Delta}
\newcommand{\La}{\Lambda}
\newcommand{\Ga}{\Gamma}

\newcommand{\LA}{\langle}
\newcommand{\RA}{\rangle}
\newcommand{\ov}[1]{\overline{#1}}
\newcommand{\un}[1]{{\underline{#1}} }

\newcommand{\Ind}{\mathop{\rm Ind}}
\newcommand{\Char}{\mathop{\rm char}}

\newcommand{\sign}{\mathop{\rm{sgn }}}
\newcommand{\Hom}{{\mathop{\rm{Hom }}}}

\newcommand{\Pf}{{\mathop{\rm{pf }}}}
\renewcommand{\P}{{\mathop{\ov{\rm{pf}}}}}
\newcommand{\DP}{{\rm DP} }
\newcommand{\F}{{\mathop{\rm{bpf }}}}

\newcommand{\NN}{{\mathbb{N}} }
\newcommand{\ZZ}{{\mathbb{Z}} }
\newcommand{\QQ}{{\mathbb{Q}} }
\newcommand{\M}{{\mathcal M} }
\newcommand{\ovphi}[1]{\varphi(#1)}

\newcommand{\QS}{\mathfrak{Q}}

\newcommand{\Q}{\mathcal{Q}}    
\newcommand{\n}{\boldsymbol{n}} 
\newcommand{\g}{\boldsymbol{g}} 
\newcommand{\h}{\boldsymbol{h}} 
\renewcommand{\i}{\boldsymbol{i}} 

\newcommand{\D}{{D}}
\newcommand{\R}{{R}}
\newcommand{\New}{0}
\renewcommand{\L}{{L}}
\newcommand{\Sp}{S\!p}

\newcommand{\Symmgr}{\mathcal{S}}                 

\begin{document}
\title{Invariants of quivers under the action of classical groups.}
 \author{
A.A. Lopatin \\
Institute of Mathematics, \\
Siberian Branch of \\
the Russian Academy of Sciences, \\
Pevtsova street, 13,\\
Omsk 644099 Russia \\
artem\underline{ }lopatin@yahoo.com \\
http://www.iitam.omsk.net.ru/\~{}lopatin/\\
}
\maketitle

\begin{abstract}
We consider a generalization of representations of quivers that
can be derived from the ordinary representations of quivers by
considering a product of arbitrary classical groups instead of a
product of the general linear groups and by considering the dual
action of groups on ``vertex"{} vector spaces together with the
usual action. A generating system for the corresponding algebra of
invariants is found. In particular, a generating system for the
algebra of $SO(n)$-invariants of several matrices is constructed
over a field of characteristic different from $2$. The proof uses
the reduction to semi-invariants of mixed representations
of a quiver and the decomposition formula that generalizes
Amitsur's formula for the determinant.
\end{abstract}



\section{Introduction}\label{section_intro}
We work over an infinite field $K$ of arbitrary characteristic.
All vector spaces, algebras, and modules are over $K$ unless
otherwise stated.

A {\it quiver} is a finite oriented graph. This notion was
introduced by Gabriel in~\cite{Gabriel72} as an effective mean
for description of different problems of the linear algebra. The
importance of this notion from the point of view of
the representational theory is due to the fact that the category of
representations of a quiver is equivalent to the category of
finite dimensional modules over the path algebra associated with
the quiver. Since every finite dimensional basic
algebra over algebraically closed field is a factor-algebra of the
path algebra of some quiver (see Chapter~3
from~\cite{Kirichenko}), the category of finite dimensional modules
over such an algebra is a full subcategory of the category of
representations of the quiver. Invariants of quivers are important
not only in the invariant theory but also in the representational
theory because these invariants distinguish semi-simple
representations of a quiver.

A representation of a quiver with $l$ vertices consists of a
collection of column vector spaces $K^{\n_1},\ldots,K^{\n_l}$,
assigned to the vertices, and linear mappings between the vector
spaces ``along"{} the arrows. We generalize this notion as
follows. Let $v$, where $1\leq v\leq l$, be a vertex of the quiver. In
the classical case $GL(\n_v)$ acts on $K^{\n_v}$ but in our
case an arbitrary classical group from the list $GL(\n_v)$,
$O(\n_v)$, $Sp(\n_v)$, $SL(\n_v)$, $SO(\n_v)$ can act on
$K^{\n_v}$. Moreover, we consider the dual
space $(K^{\n_v})^{\ast}$ together with $K^{\n_v}$
in order to deal with bilinear forms
together with linear mappings. Finally, instead of arbitrary
linear mappings ``along"{} arrows we consider only
those that, for example, preserve some bilinear symmetric form
on ``vertex"{} spaces, etc. This construction is called a {\it
mixed quiver setting}. The exact definition together with examples is given in
Section~\ref{subsection_definitions}.

{\it Orthogonal} and {\it symplectic} representations of {\it
symmetric} quivers, {\it (super)mixed} representations of quivers, and representations of {\it
signed} quivers, respectively,
introduced by Derksen and Weyman in~\cite{DW02}, Zubkov in~\cite{ZubkovI}, and
Shmelkin in~\cite{Shmelkin}, respectively,
are partial cases of this construction (see part~3 of Example~\ref{ex3}). The motivation for
these generalizations of quivers from the point of view of the representational
theory of algebraic groups was given in~\cite{DW02},~\cite{Shmelkin},
where symmetric and signed quivers, respectively, of tame and
finite type were classified.

In the paper we established generators for the invariants of a
mixed quiver setting (see Theorem~\ref{theo_main}). In particular,
in Section~\ref{subsection_matrices} we completed description,
originated by Sibirskii in~\cite{Sibirskii68} and Procesi in~\cite{Procesi76}, of generators for
the invariants of several matrices under the diagonal action by
conjugation of a classical group.

The paper is organized as follows.

Section~\ref{section_quivers} is started with the definition of
generalized quivers and their representations. It is followed by an
overview of known results on generating systems for the invariants
of quivers.

Section~\ref{section_prelim} contains notations that are
used throughout the paper. We also
recall some definitions from~\cite{Lop_bplp} such as a block
partial linearization of the pfaffian (b.p.l.p.) and a tableau with substitution.

In Section~\ref{section_geneartors} our main result
(Theorem~\ref{theo_main}) is formulated in terms of b.p.l.p.-s,
where non-zero blocks are ``generic"{} matrices. At the end of the section we consider
the quiver with one vertex
(Corollary~\ref{cor_matrices}). This special case is of great
importance since it is the simplest one and, roughly speaking, the general case can be reduced to it.

Sections~\ref{section_action}--\ref{section_proof_main} are
dedicated to the proof of Theorem~\ref{theo_main}. In
Section~\ref{section_action} we show that elements from
Theorem~\ref{theo_main} are invariants (Lemma~\ref{lemma_they are
invariants}). In Section~\ref{section_GL_SL} we rewrite generators
of semi-invariants from~\cite{LZ1} in terms of b.p.l.p.-s, and
thus obtain a set generating the algebra of $GL$- and
$SL$-invariants as a vector space over $K$
(Theorem~\ref{theo_GL_SL}). The general case is reduced to $GL$-
and $SL$-invariants in Section~\ref{section_reduction} by means of
Frobenius reciprocity and the theory of modules with good
filtration. In Theorem~\ref{theo_surjection} we show that the space of invariants is the image of
$GL$- and $SL$-invariants of an explicitly constructed quiver. In
Section~\ref{section_tableau_pairs} we apply the decomposition
formula from~\cite{Lop_bplp} to rewrite a b.p.l.p. as a polynomial
in b.p.l.p.-s of a special form
(Theorem~\ref{theo_decomposition}). Using this, in
Section~\ref{section_proof_main} we describe the image of $GL$-
and $SL$-invariants and show
that elements from Theorem~\ref{theo_main} generate the algebra of
invariants.

\section{Generalized representations of quivers}\label{section_quivers}

\subsection{Definitions}\label{subsection_definitions}
A {\it quiver} $\Q=(\Q_0,\Q_1)$ is a finite oriented graph, where
$\Q_0=\{1,\ldots,l\}$ is the set of vertices and $\Q_1$ is the set
of arrows. For an arrow $\al$, denote by $\al'$ its head and by
$\al''$ its tail. Given a {\it dimension vector}
$\n=(\n_1,\ldots,\n_l)$, we assign an $\n_v$-dimensional vector space
$V_v$ to $v\in \Q_0$. We identify $V_v$ with the space of column
vectors $K^{\n_v}$. Fix the {\it standard} basis
$e(v,1),\ldots,e(v,\n_v)$ for $K^{\n_v}$, where $e(v,i)$ is a
column vector whose $i$-th entry is $1$ and the rest of entries
are zero. A {\it representation} of $\Q$ of dimension vector
$\n$ is a collection of matrices %
$$h=(h_\al)_{\al\in \Q_1}\in %
H=H(\Q,\n)=\bigoplus_{\al\in \Q_1} K^{\n_{\al'}\times \n_{\al''}} \simeq%
\bigoplus_{\al\in \Q_1} \Hom_K(V_{\al''},V_{\al'}),$$ %
where $K^{n_1\times n_2}$ stands for the linear space of
$n_1\times n_2$ matrices over $K$ and the isomorphism is given by
the choice of bases. We will refer to $H$ as the {\it space of
representations} of $\Q$ of dimension vector $\n$.
The action of the group %
$$G=GL(\n)=\prod_{v\in \Q_0} GL(\n_v)$$ %
on $H$ is via the change of the bases for $V_v$ ($v\in \Q_0$). In
other words, $GL(\n_v)$ acts on $V_v$ by left multiplication, and
this action induces the action of $G$ on $H$ by
$$g\cdot h=(g_{\al'}h_\al g_{\al''}^{-1})_{\al\in \Q_1},$$ %
where $g=(g_\al)_{\al\in \Q_1}\in G$ and $h=(h_\al)_{\al\in \Q_1}\in H$.

The coordinate ring of the affine variety $H$ is the polynomial ring  %
$$K[H]=K[x_{ij}^\al\,|\,\al\in \Q_1,\,1\leq i\leq \n_{\al'},1\leq j\leq\n_{\al''}].$$ %
Here $x_{ij}^\al$ stands for the coordinate function on $H$ that
takes a representation $h\in H$ to the $(i,j)$-th entry of a
matrix $h_{\al}$ from $K^{\n_{\al'}\times \n_{\al''}}$. Denote by
$X_\al=(x_{ij}^\al)$ the $\n_{\al'}\times \n_{\al''}$ {\it
generic} matrix.

We will use the following notation to define the action of $G$ on
$K[H]$. Given $g\in G$, we write $g\cdot X_\al$ for the matrix,
whose $(i,j)$-th entry is $g\cdot x_{ij}^\al$. Similarly,
$\Phi(X_\al)$ stands for the matrix, whose $(i,j)$-th entry is
$\Phi(x_{ij}^\al)$, where $\Phi$ is a mapping, defined on $K[H]$.

The action of $G$ on $H$ induces the action on $K[H]$ as follows:
$(g\cdot f)(h)=f(g^{-1}\cdot h)$ for all $g\in G$, $f\in K[H]$,
$h\in H$. In other words, %
$$g\cdot X_\al=g_{\al'}^{-1}X_\al g_{\al''}$$ %
for $g\in G$, $\al\in \Q_1$. The algebra of {\it invariants} is
$$K[H]^{G}=\{f\in K[H]\,|\,g\cdot f=f\;{\rm for\; all}\;g\in G\}.$$

Given a positive integer $n$, let us fix the following notations
for the classical groups:
\begin{enumerate}\item[]
$O(n)=\{A\in K^{n\times n}\,|\,A A^t=A^t A=E\}$, $\Sp(2n)=\{A\in
K^{2n\times 2n}\,|\,A^t
J A=J\}$, $SO(n)=\{A\in O(n)\,|\,\det(A)=1\}$, where
$E=E(n)$ is the identity matrix, %
$J=J(2n)=\left(
\begin{array}{cc}
0& E(n) \\
-E(n)& 0\\
\end{array}
\right)$ is the matrix of the skew-symmetric bilinear form on
$K^{2n}$;
\end{enumerate}
and for certain subspaces of $K^{n\times n}$:
\begin{enumerate}\item[]
$S^{+}(n)=\{A\in K^{n\times n} \,|\, A^t=A\}$ is the space of symmetric matrices, %
$S^{-}(n)=\{A\in K^{n\times n} \,|\, A^t=-A\}$ is the space of
skew-symmetric matrices, %
$L^{+}(n)= %
\{A\in K^{n\times n} \,|\, AJ\text{ is a symmetric matrix}\}$, %
$L^{-}(n)=\{A\in K^{n\times n} \,|\, AJ$ is a skew-symmetric
matrix$\}$.
\end{enumerate}

The notion of representations of quivers can be generalized by the
successive realization of the following steps.
\begin{enumerate}
\item[1.] Instead of $GL(\n)$ we can take a product
$G(\n,\g)$ of classical linear groups. Here
$\g=(\g_1,\ldots,\g_l)$ is a vector, whose entries
$\g_1,\ldots,\g_l$ are symbols from the list $GL,O,\Sp,SL,SO$. By
definition,
$$G(\n,\g)=\prod_{v\in \Q_0} G_v,$$
where
$$
G_v=\left\{%
\begin{array}{ll}
GL(\n_v),& \text{if }\g_v=GL\\
O(\n_v),& \text{if }\g_v=O\\
\Sp(\n_v),& \text{if }\g_v=\Sp\\
SL(\n_v),& \text{if }\g_v=SL\\
SO(\n_v),& \text{if }\g_v=SO\\
\end{array}
\right..
$$

Obviously, we have to assume that $\n$ and $\g$ are subject to the following
restrictions:
\begin{enumerate}
\item[a)] if $\g_v=\Sp$ ($v\in \Q_0$), then $\n_v$ is even;

\item[b)] if $\g_v$ is $O$ or $SO$ ($v\in \Q_0$), then the
characteristic of $K$ is not $2$.
\end{enumerate} %

\item[2.] We can change the definition of $G(\n,\g)$ in such a
manner that allows us to deal with bilinear forms together with
linear mappings. Since bilinear forms on some vector space $V$ are
in one to one correspondence with linear mappings from the dual vector space $V^{\ast}$ to $V$,
we should change vector spaces
assigned to some vertices to the dual ones. In order to do this
consider a mapping $\i:\Q_0\to \Q_0$ such that
\begin{enumerate}
\item[c)] $\i$ is an involution, i.e., $\i^2$ is the identical
mapping;

\item[d)] $\n_{\i(v)}=\n_{v}$ for every vertex $v\in \Q_0$.
\end{enumerate}
For every $v\in \Q_0$ with $v<\i(v)$ assume that
$V_{\i(v)}=V_v^{\ast}$. Consider the dual basis
$e(v,1)^{\ast},\ldots,e(v,\n_v)^{\ast}$ for $V_v^{\ast}$ and
identify $V_v^{\ast}$ with the space of column vectors of length
$\n_v$, so $e(v,i)^{\ast}$ is the same column vector as $e(v,i)$.

The action of $GL(\n_v)$ on $V_v$ induces the action on
$V_v^{\ast}$, which we consider as the degree one homogeneous
component of the graded algebra $K[V_v]$. Given $g_v\in GL(\n_v)$
and $u\in V_v^{\ast}$, we have %
$$g_v\cdot u=(g_v^{-1})^t u.$$
Hence, we should change the group $G$ to
$$G(\n,\g,\i)=\{g\in G(\n,\g)\,|\,
g_{\i(v)}=(g_v^{-1})^t\; {\rm for\;all}\; v\in \Q_0\; {\rm
with}\;v<\i(v)\}.
$$
Since the vector spaces $K^n$ and $(K^{n})^{\ast}$ are isomorphic
as modules over $O(n)$, $\Sp(n)$, and $SO(n)$, we assume that
\begin{enumerate}
\item[e)] if $\g_v$ is $O$, $\Sp$ or $SO$ ($v\in \Q_0$), then
$\i(v)=v$.
\end{enumerate}

\item[3.] Instead of the space $H(\Q,\n)$ we should take its subspace
$H(\Q,\n,\h)$, where $\h=(\h_\al)_{\al\in \Q_1}$ and $\h_\al$ is a
symbol from the list
$M,S^{+},S^{-},L^{+},L^{-}$. By definition,  %
$$H(\Q,\n,\h)=\bigoplus_{\al\in \Q_1} H_{\al},$$
where
$$
H_{\al}=\left\{%
\begin{array}{ll}
K^{\n_{\al'}\times\n_{\al''}},& \text{if }\h_{\al}=M\\
S^{+}(\n_{\al'}),& \text{if }\h_{\al}=S^{+}\\
S^{-}(\n_{\al'}),& \text{if }\h_{\al}=S^{-}\\
L^{+}(\n_{\al'}),& \text{if }\h_{\al}=L^{+}\\
L^{-}(\n_{\al'}),& \text{if }\h_{\al}=L^{-}\\
\end{array}
\right..
$$
Additionally, we have to assume that $\n$ and $\h$ are subject to the
restriction: %
\begin{enumerate}
\item[f)] if $\h_\al\neq M$ ($\al\in \Q_1$), then
$\n_{\al'}=\n_{\al''}$.
\end{enumerate}
\end{enumerate}
Consider a group $G=G(\n,\g,\i)\subset G(\n)$ and a vector space
$H=H(\Q,\n,\h)\subset H(\Q,\n)$ satisfying the previous conditions~a)--f). To ensure that these inclusions
induce the action of $G$ on $H$, we assume that the following additional conditions are also valid for all
$v\in \Q_0$, $\al\in \Q_1$:
\begin{enumerate}\item[]
\begin{enumerate}
\item[g)] if $\al$ is a loop, i.e., $\al'=\al''$, and $\h_\al$ is
$S^{+}$ or $S^{-}$, then $\g_{\al'}$ is $O$ or $SO$;

\item[h)] if $\al$ is a loop and $\h_\al$ is $L^{+}$ or $L^{-}$,
then $\g_{\al'}=\Sp$;

\item[i)] if $\al$ is not a loop and $\h_\al\neq M$, then
$\i(\al')=\al''$ and $\h_\al$ is $S^{+}$ or $S^{-}$.
\end{enumerate}
\end{enumerate}
A quintuple $\QS=(\Q,\n,\g,\h,\i)$ satisfying~a)--i) is called a
{\it mixed quiver setting}.
Definitions of the generic
matrices $X_\al$ and the algebra of invariants $K[H]^G$ are the
same as above. Note
that %
if $\h_\al=S^{+}$, then $X_{\al}^t=X_\al$; %
if $\h_\al=S^{-}$, then $X_{\al}^t=-X_\al$; %
if $\h_\al=L^{+}$, then $(X_{\al}J)^t=X_\al J$; and %
if $\h_\al=L^{-}$, then $(X_{\al}J)^t=-X_\al J$. %
In this paper we establish a generating system for $K[H]^G$.

\example\label{ex3} {\bf 1.} Let $\Q$ be the following quiver
$$\vcenter
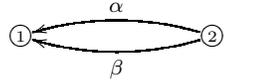
{
\xymatrix@C=1cm@R=1cm{ %
\vtx{1}\ar@/^/@{<-}[rr]^{\al} \ar@/_/@{<-}[rr]_{\be}&&\vtx{2}\\
}} \quad.
$$
Define a mixed quiver setting $\QS=(\Q,\n,\g,\h,\i)$ by $\n_1=\n_2=n$, $\g_1=\g_2=GL$,
$\h_\al=\h_\be=S^{+}$, and $\i(1)=2$. The group $G(\n,\g,\i)\simeq GL(n)$ acts on
$H(\Q,\n,\h)=S^{+}(n)\oplus S^{+}(n)$ by the rule
$$g\cdot (A,B)=(g A g^t, g B g^t)$$
for $g\in GL(n)$ and $(A,B)\in S^{+}(n)\oplus S^{+}(n)$. Hence
the orbits of this action correspond to pairs
of symmetric bilinear forms on $K^{n}$. If we put $\h_\al=\h_\be=S^{-}$, then we obtain
pairs of skew-symmetric bilinear forms on $K^{n}$.  The classification problem for such
pairs is a classical topic going back to Weierstrass and Kronecker
(see~\cite{Hodge_Pedoe_I},~\cite{Hodge_Pedoe_II}, and~\cite{Gantmacher}).
\smallskip

{\bf 2.} Let $\Q$ be the following quiver
$$\vcenter
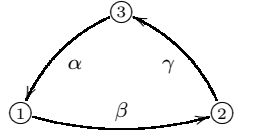
{
\xymatrix@C=1cm@R=1cm{ %
&\vtx{3}\ar@/_/@{->}[dl]^{\al} \ar@/^/@{<-}[dr]_{\ga}&\\
\vtx{1}\ar@/_/@{->}[rr]^{\be} &&\vtx{2}\\
}} \quad.
$$
Define a mixed quiver setting $\QS=(\Q,\n,\g,\h,\i)$ by $\n_1=\n_2=n$, $\n_3=m$; $\g_1=\g_2=SL$, $\g_3=O$;
$\h_\al=\h_\ga=M$, $\h_\be=S^{+}$; and $\i(1)=2$, $\i(3)=3$. Hence the action of $G(\n,\g,\i)\simeq SL(n)\times O(m)$ on
$H(\Q,\n,\h)=K^{n\times m}\oplus S^{+}(n)\oplus K^{m\times n}$ is given by
$$(g,f)\cdot (A,B,C)=(gAf^t, (g^{-1})^t B g^{-1}, f C g^t)$$
for $(g,f)\in SL(n)\times O(m)$ and $(A,B,C)\in K^{n\times m}\oplus S^{+}(n)\oplus K^{m\times n}$.
\smallskip

{\bf 3.} If we consider mixed quiver settings $\QS=(\Q,\n,\g,\h,\i)$ with the restriction
$\g_v\in\{GL,O,\Sp\}$ for all $v\in\Q_0$, then we obtain the definition of supermixed
representations of a quiver (see~\cite{ZubkovI}), or, equivalently, the definition of
representations of a signed quiver (see~\cite{Shmelkin}).

\subsection{Known results}\label{subsection_history}

In this section $(\Q,\n,\g,\h,\i)$ is a mixed quiver setting,
$G=G(\n,\g,\i)$, and $H=H(\Q,\n,\h)$. We overview the known
results on generators and relations between them for the algebra
of invariants $K[H]^G$.

The first results on invariants of quivers (i.e. for the identical involution $\i$) were
obtained for an important special case of a quiver with one vertex and several loops. For
a field of characteristic zero generators for $G\in\{GL(n),SL(n),\Sp(n)\}$
were described by Sibirskii
in~\cite{Sibirskii68} and Procesi in~\cite{Procesi76}.
Procesi also described relations between generators in~\cite{Procesi76}
applying the classical theory of invariants of vectors and
covectors (see book~\cite{Weyl46} by Weyl). Independently, relations for
$G=GL(n)$ were described by Razmyslov in~\cite{Razmyslov74}. Developing ideas
from~\cite{Procesi76}, Aslaksen et al. calculated generators for
$G=SO(n)$ (see~\cite{Aslaksen95}).

Invariants for an arbitrary quiver $\Q$ were considered by Le
Bruyn and Procesi in~\cite{Le_Bruyn_Procesi_90}, where generators
were found for the case of $G=\prod_{v\in\Q_0}GL(\n_n)$ and
$H=\sum_{\al\in\Q_1}K^{n_{\al'}\times n_{\al''}}$. Similar
results were obtained later by Domokos in~\cite{Domokos98}.

The importance of characteristic-free approach to quiver invariants
was pointed out by Formanek in overview~\cite{Formanek91} (see also~\cite{Formanek87}).
Relying on the
theory of modules with good filtrations (see~\cite{Donkin85}), Donkin
described generators for a quiver with one vertex in~\cite{Donkin92a}
and for an arbitrary quiver afterwards (see~\cite{Donkin94}). Relations
between generators from the mentioned papers were found by Zubkov
in~\cite{Zubkov96} and~\cite{Zubkov_Fund_Math_01} by means of an
approach that allowed to calculate generators and relations
between them simultaneously. His method is also based on the
theory of modules with good filtrations.

For the rest of classical groups over a field of positive
characteristic, the first results were obtained by Zubkov.
In~\cite{Zubkov99} he obtained generators for a quiver with one vertex and
the orthogonal or symplectic group $G$. The proof is based on ideas from~\cite{Donkin92a} and a
reduction to invariants of mixed quiver settings with $\g_v=GL$
for every vertex $v$. The reduction was performed by means
Frobenius reciprocity. Let us recall that we do not consider the
case of the
(special) orthogonal group in characteristic $2$ case. The reason is that in the later case
even generators of invariants of several vectors are not known
(for the latest developments see~\cite{Domokos_Frenkel}).

Invariants of a quiver under the action of
$G=\prod_{v\in\Q_0}SL(n_v)$ are called semi-invariants. Its
generators for an arbitrary characteristic were established by
Domokos and Zubkov in~\cite{DZ01} using the methods
from~\cite{Donkin92a},~\cite{Donkin94},~\cite{Zubkov96},~\cite{Zubkov_Fund_Math_01},
and, independently, by Derksen and
Weyman in~\cite{DW_LR_02},~\cite{DW00} utilizing the methods of the
representation theory of quivers. Simultaneously, similar result
in the case of characteristic zero was obtained by Schofield and
Van den Bergh in~\cite{Schofield_Van_den_Bergh_01}. These results
were generalized for mixed quiver settings with $\g_v=SL$ and
$\h_{\al}=M$, where $v$ is a vertex and $\al$ is an arrow, by
the author and Zubkov in~\cite{LZ1}.

Zubkov in~\cite{ZubkovI} combined the arguments for Young superclasses
from~\cite{Donkin94} with the reduction from~\cite{Zubkov99} to describe
generators for a mixed quiver setting $(\Q,\n,\g,\h,\i)$
with $\g_v\in\{GL,O,\Sp\}$ for all $v\in\Q_0$. Relations between
them were found in~\cite{ZubkovI}.

Generators for the remaining mixed quiver settings are computed in
this paper.

\section{Preliminaries}\label{section_prelim}
\subsection{Notations}\label{subsection_notations} %
In what follows, ${\NN}$ stands for the set of non-negative integers, $\ZZ$ for the set of integers,
and $\QQ$ for the quotient field of the ring $\ZZ$.

The cardinality of a set $S$ is denoted by $\#S$ and the permutation group on $n$ elements
is denoted by $\Symmgr_n$. Given integers $i<j$, we write $[i,j]$ for the interval $i,i+1,\ldots,j-1,j$.

By a {\it distribution} $B=(B_1,\ldots,B_s)$ of a set $[1,t]$ we
mean an ordered partition of the set into pairwise disjoint
subsets $B_i$ ($1\leq i\leq s$), which are called components of
the distribution. To every $B$ we associate two functions
$j\mapsto B|j|$ and $j\mapsto B\LA j\RA$ ($1\leq j\leq t$),
defined by the rules:
$$B|j|=i, \text{ if }j\in B_i,\text{ and }B\LA
j\RA=\#\{[1,j]\cap B_i\,|\,j\in B_i\}.$$ %

A vector $\un{t}=(t_1,\ldots,t_s)\in \NN^s$ determines the
distribution $T=(T_1,\ldots,T_s)$ of the set $[1,t]$, where
$t=t_1+\cdots+t_s$ and
$T_i=\{t_1+\cdots+t_{i-1}+1,\ldots,t_1+\cdots+t_i\}$, $1\leq i\leq
s$. As an example, if $\un{t}=(1,3,0,2)$, then
$T=(\{1\},\{2,3,4\},\emptyset,\{5,6\})$ and $T|5|=4$.

A vector $\un{\la}=(\la_1,\ldots,\la_s)\in\NN^s$ satisfying
$\la_1\geq\cdots \geq\la_s$ and $\la_1+\cdots+\la_s=t$ is called a
{\it partition} of $t$ and is denoted by $\un{\la}\vdash t$. A
{\it multi-partition} $\un{\la}\vdash \un{t}$ is a $q$-tuple of
partitions $\un{\la}=(\un{\la}_1,\ldots,\un{\la}_q)$, where
$\un{\la}_1\vdash t_1,\ldots,\un{\la}_q\vdash t_q$, and
$\un{t}=(t_1,\ldots,t_q)\in\NN^q$.

\subsection{Pfaffians and tableaux with substitutions}\label{subsection_pf}

Denote coefficients in the characteristic polynomial
of an $n\times n$ matrix $X$ by $\si_k(X)$, i.e., %
$$\det(\la E-X)=\la^n-\si_1(X)\la^{n-1}+\cdots+(-1)^n\si_n(X).$$ %

Assume $n$ is even. Define the {\it generalized
pfaffian} of an arbitrary $n\times n$ matrix $X=(x_{ij})$ by
$$\P(X)=\Pf(X-X^{t}),$$
\noindent where $\Pf$ stands for the pfaffian of a skew-symmetric matrix. By
abuse of notation we will refer to $\P$ as the pfaffian. For $K=\QQ$ there is a more
convenient formula
\begin{eq}\label{eq_P}
\P(X)=\Pf(X-X^{t})=\frac{1}{(n/2)!}\sum\limits_{\pi\in \Symmgr_{n}}\sign(\pi)
\prod\limits_{i=1}^{n/2} x_{\pi(2i-1),\pi(2i)}.%
\end{eq} %

For $n\times n$ matrices $X_1=(x_{ij}(1)),\ldots,X_s=(x_{ij}(s))$
and positive integers $k_1,\ldots,k_s$, satisfying
$k_1+\cdots+k_s=n/2$, consider the polynomial $\P(x_1
X_1+\cdots+x_s X_s)$ in the variables $x_1,\ldots,x_s$. The
partial linearization $\P_{k_1,\ldots,k_s}(X_1,\ldots,X_s)$ of the
pfaffian is the coefficient at $x_1^{k_1}\cdots x_s^{k_s}$ in this
polynomial. Assume that for some $\un{n}=(n_1,\ldots,n_m)\in\NN^m$
with $n_1+\cdots+n_m=n$ each of matrices $X_1,\ldots,X_s$ is
partitioned into $m\times m$ number of blocks, where the block in
the $(i,j)$-th position is an $n_i\times n_j$ matrix and the only
non-zero block is the one in the $(p,q)$-th position. Then
$\P_{k_1,\ldots,k_s}(X_1,\ldots,X_s)$ is called a {\it block
partial linearization of the pfaffian} ({\it b.p.l.p.}).

The following notions were introduced in Section~3 of~\cite{Lop_bplp},
where more detailed explanation and examples are given.

\bigskip
\noindent{\bf Definition (of shapes).} The {\it shape} of
dimension $\un{n}=(n_1,\ldots,n_m)\in\NN^{m}$ is the collection of
$m$ columns of cells. The columns are numbered by $1,2,\ldots,m$,
and the $i$-th column contains exactly $n_i$ cells, where $1\leq
i\leq m$. Numbers $1,\ldots,n_i$ are assigned to the cells of the
$i$-th column, starting from the top. As an example, the shape of
dimension $\un{n}=(3,2,3,1,1)$ is%
$$\begin{tabular}{ccccc} %
1&2&3&4&5  \\
\end{tabular}$$ %
\vspace{-.7cm}
$$\begin{tabular}{|c|c|c|c|c|} %
\hline 1&1&1&1&1       \\
\hline 2&2&2         \\
\hhline{---~~} 3&&3   \\
\hhline{-~-~~}
\end{tabular}$$ %

\bigskip
\noindent{\bf Definition (of a tableau with substitution).} Let
$\un{n}=(n_1,\ldots,n_m)\in\NN^{m}$ and let $n=n_1+\cdots+n_m$ be
even. A pair $(T,(X_1,\ldots,X_s))$ is called a {\it tableau with substitution}
of dimension $\un{n}$ if
\begin{enumerate}
\item[$\bullet$] $T$ is the shape of dimension $\un{n}$ together
with a set of arrows. An {\it arrow} goes from one cell of the shape
into another one, and each cell of the shape is either the head or
the tail of one and only one arrow.  We refer to $T$ as {\it
tableau} of dimension $\un{n}$, and we write $a\in T$ for an arrow
$a$ from $T$. Given an arrow $a\in T$, denote by $a'$ and $a''$
the columns containing the head and the tail of $a$, respectively.
Similarly, denote by $'a$ the number assigned to the cell
containing the head of $a$, and denote by $''a$ the number
assigned to the cell containing the tail of $a$. Schematically
this is depicted as
$$\begin{array}{cccc} %
\;\;&a''&\! a'&  \\
\end{array}
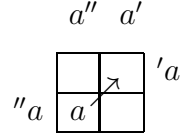
$$ %
\vspace{-.6cm}
$$\begin{array}{c|c|c|c} %
\hhline{~--~} \vspace{-0.25cm} &&&{}'a\\
              \vspace{-0.25cm} &&\hspace{-0.45cm} \nearrow&\\
\hhline{~--~}''a &a&\;\;&\\
\hhline{~--~}
\end{array}
$$%

\item[$\bullet$] $\varphi$ is a fixed mapping from the set of arrows of $T$ onto
$[1,s]$ that satisfies the following property:

\begin{enumerate}
\item[] if $a,b\in T$ and $\ovphi{a}=\ovphi{b}$,
then $a'=b'$, $a''=b''$;
\end{enumerate}

\item[$\bullet$] $(X_1,\ldots,X_s)$ is a sequence of matrices
such that the matrix $X_{\ovphi{a}}$ assigned to the arrow $a\in T$ is
$n_{a''}\times n_{a'}$ matrix and its $(p,q)$-th entry is denoted by $(X_j)_{pq}$.
\end{enumerate}

\example\label{ex_new} Let $T$ be the tableau
$$\begin{array}{|c|c|c|}
\hline \vspace{-0.2cm} a & b &\;\,\\
\vspace{-0.25cm}\downarrow &\searrow\hspace{-0.60cm}
&\hspace{-0.50cm}\nearrow \\ 
\hline \vspace{-0.25cm}& c&\\
\vspace{-0.25cm}& \\    
\hline d-\!\!\!\! &\!\!\!\!\rightarrow\\
\hhline{--~}
\end{array}
$$ %
of dimension $(3,3,2)$. Define $\varphi$ by $\ovphi{a}=1$, $\ovphi{b}=\ovphi{c}=2$, and $\ovphi{d}=3$, and let
$X_1$, $X_3$ be $3\times 3$ matrices and $X_2$ be a $3\times 2$ matrix. Then $(T,(X_1,X_2,X_3))$ is a tableau with substitution.

\bigskip
\noindent{\bf Definition (of ${\F}_T(X_1,\ldots,X_s))$.} 
Let $(T,(X_1,\ldots,X_s))$ be a tableau with substitution of dimension
$\un{n}$. Define the polynomial %
\begin{eq}\label{eq_def_F0}
{\F}_T^0(X_1,\ldots,X_s)=\sum_{\pi_1\in
\Symmgr_{n_1},\ldots,\pi_m\in \Symmgr_{n_m}}
\sign(\pi_1)\cdots\sign(\pi_m)\prod_{a\in T}
(X_{\ovphi{a}})_{\pi_{a''}(''a),\pi_{a'}('a)},
\end{eq}%
and the coefficient
$$c_T=\prod_{j=1}^s\#\{a\in T\,|\,\ovphi{a}=j\}!$$ %
In the case $K=\QQ$ define
$${\F}_T(X_1,\ldots,X_s)=\frac{1}{c_T}{\F}^0_T(X_1,\ldots,X_s).$$ %
Since ${\F}_T(X_1,\ldots,X_s)$ is a polynomial in entries of
$X_1,\ldots,X_s$ with integer coefficients, the definition of ${\F}_T(X_1,\ldots,X_s)$
extends over an arbitrary field.

\example\label{ex1} For every $n\times n$ matrix $X$ there is a
tableau with substitution $(T,X)$ such that $\det(X)={\F}_T(X)$. If $n$ is odd,
then the same is valid for $\P(X)$ (see Example~2 of~\cite{Lop_bplp}
for details).

\bigskip
\noindent{}The next lemma, which is part~b) of
Lemma~1 from~\cite{Lop_bplp}, shows that $\F$ is a b.p.l.p.
\begin{lemma}\label{lemma_bplp}
Let $(T,(X_1,\ldots,X_s))$ be a tableau with substitution of dimension
$\un{n}=(n_1,\ldots,n_m)$. Consider $a_1,\ldots,a_s\in T$ such
that $\ovphi{a_1}=1,\ldots,\ovphi{s}=s$. For any $1\leq p\leq s$
denote by $Z_p$ the $n\times n$ matrix, partitioned into $m\times
m$ number of blocks, where the block in the $(i,j)$-th position is
an $n_i\times n_j$ matrix; the block in the $(a_p'',a_p')$-th
position is equal to $X_p$, and the rest of blocks are zero
matrices. Then
$$\F_T(X_1,\ldots,X_s)=\pm\P_{k_1,\ldots,k_s}(Z_1,\ldots,Z_s),$$
where $k_p=\#\{a\in T\,|\,\ovphi{a}=p\}$ for any $1\leq p\leq s$.
\end{lemma}

%
%

\section{Generators}\label{section_geneartors}

\subsection{Main results}\label{subsection_results} %
Let $\QS=(\Q,\n,\g,\h,\i)$ be a mixed quiver setting and
$\Q_0=\{1,\ldots,l\}$. Without loss of generality we can assume
that
\begin{eq}\label{eq_condition}
{\rm\; if\;} v\in \Q_0 {\rm\; and\;} \g_v{\rm\; is\;} GL {\rm\;
or\;}SL, {\rm\; then\;} \i(v)\neq v.
\end{eq}%
Otherwise we can add a new vertex $\ov{v}$ to $\Q$, and set
$\i(v)=\ov{v}$, $\n_{\ov{v}}=\n_v$, $\g_{\ov{v}}=\g_{v}$; this
construction changes neither the space $H(\Q,\n,\h)$ nor
the algebra of invariants.

\bigskip
\noindent\textbf{ Definition (of the mixed double quiver setting
$\QS^{\D}$).} %
Define the mixed {\it double} quiver setting
$\QS^{\D}=(\Q^{\D},\n,\g,\h^{\D},\i)$ as follows:
$\Q_0^{\D}=\Q_0$, $\Q_1^{\D}=\Q_1\coprod \{\al^t\,|\,\al\in
\Q_1,\; \h_\al=M\}$, where $(\al^t)'=\i({\al''})$,
$(\al^t)''=\i({\al'})$, and $\h_{\al^t}^{\D}=M$ for $\al\in \Q_1$ with $\h_\al=M$ and
$\h_{\al}^{\D}=\h_{\al}$ for all $\al\in\Q_1$.

Define a mapping $\Phi^{\D}:K[H(\Q^{\D},\n,\h^{\D})]\to
K[H(\Q,\n,\h)]$ such that
$\Phi^{\D}(X_{\al})=X_\al$ for $\al\in \Q_1$, and $\Phi^{\D}(X_{\al^t})$ for $\al\in \Q_1$
and $\h_\al=M$ is defined as follows:
\begin{enumerate} %
\item[$\bullet$] If $\g_{\al'}\neq \Sp$ and $\g_{\al''}\neq \Sp$,
then $\Phi^{\D}(X_{\al^t})=X_\al^t$.

\item[$\bullet$] If $\g_{\al'}=\Sp$ and $\g_{\al''}\neq \Sp$, then
$\Phi^{\D}(X_{\al^t})=X_\al^t J(\n_{\al'})$.

\item[$\bullet$] If $\g_{\al'}\neq \Sp$ and $\g_{\al''}=\Sp$, then
$\Phi^{\D}(X_{\al^t})=J(\n_{\al''})X_\al^t$.

\item[$\bullet$] If $\g_{\al'}=\Sp$ and
$\g_{\al''}=\Sp$, then %
$\Phi^{\D}(X_{\al^t})=J(\n_{\al''}) X_\al^t J(\n_{\al'}) $.
\end{enumerate} %
Let us remark that the meaning of notation $\Phi^{\D}(X_{\al})$ was
explained in Section~\ref{subsection_definitions}.

\example\label{ex2} Let $\Q$ be
$$\vcenter
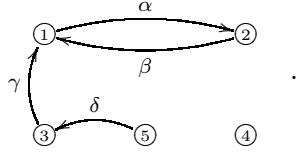
{
\xymatrix@C=1cm@R=1cm{ %
\vtx{1}\ar@/^/@{->}[rr]^{\al} \ar@/_/@{<-}[rr]_{\be}&&\vtx{2}\\
\vtx{3}\ar@/^/@{<-}[r]^{\de} \ar@/^/@{->}[u]^{\ga}&\vtx{5}& \vtx{4}\\
}} \quad.
$$
Define a mixed quiver setting $\QS=(\Q,\n,\g,\h,\i)$ by $\i(1)=2$,
$\i(3)=4$, $\i(5)=5$; $\g_1=\g_2=GL$, $\g_3=\g_4=SL$, $\g_5=O$;
$\h_\al=\h_\ga=\h_\de=M$, $\h_\be=S^{+}$. Then $\Q^{\D}$ is
$$
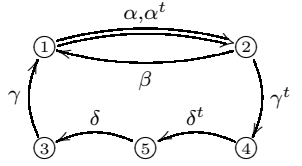
\vcenter{
\xymatrix@C=1cm@R=1cm{ %
\vtx{1}\ar@2@/^/@{->}[rr]^{\al,\al^t} \ar@/_/@{<-}[rr]_{\be}&&\vtx{2}\\
\vtx{3}\ar@/^/@{<-}[r]^{\de} \ar@/^/@{->}[u]^{\ga} &\vtx{5}&
\vtx{4}\ar@/_/@{->}[l]_{\de^t} \ar@/_/@{<-}[u]_{\ga^t} \\
}} .
$$
\smallskip

Before presenting the next concept, let us recall that
$\al=\al_1\cdots \al_r$ is a {\it path} in $\Q$ (where
$\al_1,\ldots,\al_r\in \Q_1$), if
$\al_1'=\al_2'',\ldots,\al_{r-1}'=\al_r''$. The head of the path
$\al$ is $\al'=\al_r'$ and the tail is $\al''=\al_1''$.  The path
$\al$ is called {\it closed} if $\al'=\al''$.
\bigskip

\noindent\textbf{ Definition (of a $\QS$-tableau with substitution and a path
$\QS$-tableau with substitution).} A tableau with substitution $(T,(Y_1,\ldots,Y_s))$ of
dimension $\un{n}\in\NN^m$ is called a {\it $\QS$-tableau with substitution},
if for some {\it weight} $\un{w}=(w_1,\ldots,w_l)\in\NN^l$ and the
distribution $W$, determined by $\un{w}$ (see
Section~\ref{subsection_notations}), we have
\begin{enumerate}
\item[$\bullet$]
$\un{n}=(\underbrace{\n_1,\ldots,\n_1}_{w_1},\ldots,\underbrace{\n_l\ldots,\n_l}_{w_l})$;

\item[$\bullet$] if $a\in T$, then there exists an $\al\in \Q_1$
such that $Y_{\ovphi{a}}=X_\al$, $W|a'|=\i(\al'')$, $W|a''|=\al'$.
\end{enumerate}
If we replace the last condition by the following one
\begin{enumerate}
\item[$\bullet$] if $a\in T$, then there exists a path
$\al=\al_1\cdots \al_r$ in $\Q$ (where $\al_1,\ldots,\al_r\in
\Q_1$) such that $Y_{\ovphi{a}}=X_{\al_r}\cdots X_{\al_1}$,
$W|a'|=\i(\al'')$, $W|a''|=\al'$,
\end{enumerate}
then we obtain the definition of a {\it path $\QS$-tableau with substitution}. Obviously,
for a (path) $\QS$-tableau with substitution $(T,(Y_1,\ldots,Y_s))$ we have
$\F_T(Y_1,\ldots,Y_s)\in K[H(\Q,\n,\h)]$.

\begin{theo}\label{theo_main} %
\textbf{ (Main theorem)} Let $(\Q,\n,\g,\h,\i)$ be a mixed quiver
setting satisfying~\Ref{eq_condition}. Then the algebra of
invariants $K[H(\Q,\n,\h)]^{G(\n,\g,\i)}$ is generated as
$K$-algebra by the elements $\Phi^{\D}(\si_k(X_{\be_r}\cdots
X_{\be_1}))$, $\Phi^{\D}(\F_T(Y_1,\ldots,Y_s))$, where
\begin{enumerate}
\item[1.] $\be_1\cdots \be_r$ ranges over all closed paths in $\Q^{\D}$ and
$1\leq k\leq \n_{\be''_1}$;    

\item[2.] $(T,(Y_1,\ldots,Y_s))$ ranges over all path $\QS^{\D}$-tableaux with substitutions
of a weight $\un{w}$ such that 

\begin{enumerate} %
\item[a)] if $\g_v\in\{GL,O,\Sp\}$ for some $v\in \Q_0$, then
$w_{\i(v)}=w_{v}=0$;  

\item[b)] if $\g_v=SL$ for some $v\in \Q_0$, then $w_{\i(v)}=0$ or
$w_{v}=0$;            

\item[c)] if $\g_v=SO$ for some $v\in \Q_0$, then $w_{v}\leq 1$ and
$\i(v)=v$.            
\end{enumerate}
\end{enumerate}
\end{theo}

This theorem implies the main result of~\cite{ZubkovI}.
\begin{cor}\label{cor_ZubI}
Let $(\Q,\n,\g,\h,\i)$ be a mixed quiver setting satisfying~\Ref{eq_condition}.
If $\g_v\in\{GL,O,\Sp\}$ for all $v\in
\Q_0$, then $K$-algebra $K[H(\Q,\n,\h)]^{G(\n,\g,\i)}$ is
generated by $\Phi^{\D}(\si_k(X_{\be_r}\cdots X_{\be_1}))$, where
$\be_1\cdots \be_r$ is a closed path in $\Q^{\D}$ and $1\leq k\leq
\n_{\be''_1}$.
\end{cor}

For $f\in K[H(\Q,\n,\h)]$ let $\un{t}=(t_{\al})_{\al\in
\Q_1}\in\NN^{\#\Q_1}$ be the {\it multidegree} of $f$, i.e.,
$t_{\al}$ is the total degree of the polynomial $f$ in
$x_{ij}^{\al}$, where $1\leq i\leq \n_{\al'}$ and $1\leq j\leq
\n_{\al''}$. The algebra of invariants is homogeneous with respect to the grading by multidegrees
as well as the generating system from Theorem~\ref{theo_main}.

\subsection{Invariants of several matrices}\label{subsection_matrices}
Consider the case of a quiver with one vertex and $d$ loops. Let
$\h_\al=M$ for every arrow $\al$. Then $H=K^{n\times
n}\oplus\cdots\oplus K^{n\times n}$ is $d$-tuple of $n\times n$
matrices over $K$, and $G$ is a group from the list $GL(n)$,
$O(n)$, $\Sp(n)$, $SO(n)$. We assume that if $G$ is $O(n)$ or
$SO(n)$, then the characteristic of $K$ is nor $2$; if $G$ is
$\Sp(n)$, then $n$ is even. The group $G$ acts on $H$ by the
diagonal conjugation. Hence it acts on $K[H]$ as follows: $g\cdot
X_{\al}=g^{-1}X_\al g$, where $g\in G$, $1\leq\al\leq d$, and
$X_\al=(x_{ij}^\al)$ is the $n\times n$ generic matrix. We do not
consider the case $G=SL(n)$ because invariants for
$GL(n)$ and $SL(n)$ are the same.

\begin{cor}\label{cor_matrices}
The algebra of invariants $K[H]^G$ is generated by the following
elements:
\begin{enumerate}
\item[a)] $\si_k(X)$ ($1\leq k\leq n$ and $X$ ranges over all monomials in $X_1,\ldots,X_d$), if $G=GL(n)$;

\item[b)] $\si_k(Y)$ ($1\leq k\leq n$), if $G=O(n)$;

\item[c)] $\si_k(Y)$ ($1\leq k\leq n$), if $G=SO(n)$ and $n$ is
odd;

\item[d)] $\si_k(Y)$, $\P_{k_1,\ldots,k_s}(Y_1,\ldots,Y_s)$
($1\leq k\leq n$, $k_1+\cdots+k_s=n/2$),  if $G=SO(n)$ and $n$ is
even.
\end{enumerate}

In~b),~c), and~d) matrices $Y,Y_1,\ldots,Y_s$ range over all
monomials in $X_1,\ldots,X_d$, $X_1^t,\ldots,X_d^t$.

\begin{enumerate}
\item[e)] $\si_k(Z)$ ($1\leq k\leq n$ and $Z$ ranges over all
monomials in $X_1,\ldots,X_d$, $JX_1^tJ,\ldots,JX_d^tJ$), if $G=\Sp(n)$.
\end{enumerate}
\end{cor}
\begin{proof} It follows immediately from
Theorem~\ref{theo_main}. To prove part~d) we should also use the
fact that $\F_T(Y_1,\ldots,Y_s)$ is a partial linearization of the
pfaffian for one column tableau $T$ (see part~1 of Example~2
from~\cite{Lop_bplp}).
\end{proof}

The only new part of Corollary~\ref{cor_matrices} is part~d) (see
Section~\ref{subsection_history} for references).

\section{Action of groups on tableaux with substitutions}\label{section_action}
In this section we show that the elements from
Theorem~\ref{theo_main} are invariants.

Consider an $\un{n}\in\NN^{m}$ and the group $G=\prod_{i=1}^m
GL(n_i)$. Given $g=(g_i)_{1\leq i\leq m}\in G$ and a tableau with substitution
$(T,(Y_1,\ldots,Y_s))$ of dimension $\un{n}$, define the tableau with substitution $(T,(g* Y_1,\ldots,g* Y_s))$ by $g* Y_j=g_{a''} Y_j
g_{a'}^t$, where $1\leq j\leq s$ and $a\in T$ is such that $\ovphi{a}=j$.
Since for all $a,b\in T$ with $\ovphi{a}=\ovphi{b}$ we have
$a'=b'$ and $a''=b''$, the matrix $g* Y_j$ is well defined.

\begin{lemma}\label{lemma_GactsonT}
Using the preceding notation we have the equality
$$\F_T(g* Y_1,\ldots,g* Y_s)=\det(g_1)\cdots\det(g_m)\F_T(Y_1,\ldots,Y_s).$$
\end{lemma}
\begin{proof}
Let $n=n_1+\cdots+n_m$ and let $g_0$ be an $n\times n$ block-diagonal matrix such that
the $i$-th block is equal to $g_i$ ($1\leq i\leq m$). By the definition of $T$, $n$ is
even and $\{\ovphi{a}\,|\,a\in T\}=[1,s]$.

Repeat construction from Lemma~\ref{lemma_bplp}. Consider
$a_1,\ldots,a_s\in T$ such that
$\ovphi{a_1}=1,\ldots,\ovphi{s}=s$. For any $1\leq p\leq s$ denote
by $Z_p$ the $n\times n$ matrix, partitioned into $m\times m$
number of blocks, where the block in the $(i,j)$-th position is an
$n_i\times n_j$ matrix; the block in the $(a_p'',a_p')$-th
position is equal to $Y_p$, and the rest of blocks are zero
matrices. Then
$\F_T(Y_1,\ldots,Y_s)=q\P_{k_1,\ldots,k_s}(Z_1,\ldots,Z_s)$, where
$q=\pm1$, $k_p=\#\{a\in T\,|\,\ovphi{a}=p\}$ for any $1\leq
p\leq s$. Thus, %
$$\F_T(g* Y_1,\ldots,g* Y_s)=
q\P_{k_1,\ldots,k_s}(g_0 Z_1 g_0^t,\ldots,g_0 Z_s g_0^t)$$
$$=q\det(g_0)\P_{k_1,\ldots,k_s}(Z_1,\ldots,Z_s),$$ %
since $\P(g_0 Z g_0^t)=\det(g_0)\P(Z)$ for every $n\times n$
matrix $Z$. This completes the proof.
\end{proof}

\begin{lemma}\label{lemma_action}
Let $(\Q,\n,\g,\h,\i)$ be a mixed quiver setting that
satisfies~\Ref{eq_condition} and $\Q_0=\{1,\ldots,l\}$. Let
$(T,(Y_1,\ldots,Y_s))$ be a $\QS$-tableau with substitution or path
$\QS$-tableau with substitution of a weight $\un{w}=(w_1,\ldots,w_l)$. Assume
that $w_v=0$ for all $v\in \Q_0$ with $\g_v=\Sp$. Then for all
$g=(g_v)_{v\in \Q_0}\in G(\n,\g,\i)$ we have
$$g\cdot \F_T(Y_1,\ldots,Y_s)=q\, \F_T(Y_1,\ldots,Y_s),$$
where %
$$q=
\prod_{v\in \Q_0,\,v=\i(v)}\det(g_v)^{-w_v}\; %
\prod_{v\in \Q_0,\,v<\i(v)}\det(g_v)^{w_{\i(v)}-w_v}.
$$ %
\end{lemma}
\begin{proof} {\bf a)} Assume that $(T,(Y_1,\ldots,Y_s))$ is a $\QS$-tableau with substitution.
Let $w=w_1+\cdots+w_l$ and let $W$ be the distribution determined
by $\un{w}$. Define a mapping $\psi:G(\n,\g,\i)\to\prod_{j=1}^w
GL(\n_{W|j|})$ by $\psi(g)=(\psi(g)_1,\ldots,\psi(g)_w)$ and
$\psi(g)_j=g_{W|j|}$ for $1\leq j\leq w$. We claim
that
\begin{eq}\label{eq_gF}
g\cdot \F_T(Y_1,\ldots,Y_s)=\F_T(\psi(g^{-1})*
Y_1,\ldots,\psi(g^{-1})* Y_s) \text{ for any }g\in G(\n,\g,\i).
\end{eq}
\noindent By the definition of a $\QS$-tableau with substitution, for every $a\in T$
there exists an arrow $\al\in \Q_1$ with $Y_{\ovphi{a}}=X_\al$, $W|a'|=\i(\al'')$, and
$W|a''|=\al'$. Hence $g\cdot Y_{\ovphi{a}}=g\cdot X_\al=g_{\al'}^{-1}X_\al g_{\al''}$. On
the other hand,
$\psi(g^{-1})*Y_{\ovphi{a}}=\psi(g^{-1})_{a''}X_\al\psi(g^{-1})_{a'}^t=%
g_{W|a''|}^{-1}X_\al (g_{W|a'|}^{-1})^t=%
g_{\al'}^{-1}X_\al (g_{\i(\al'')}^{-1})^t=g_{\al'}^{-1}X_\al
g_{\al''}$, since $g_{\i(\al'')}=(g_{\al''}^{-1})^t$ for
$\g_{\al''}\neq \Sp$. Therefore $g\cdot
Y_{\ovphi{a}}=\psi(g^{-1})*Y_{\ovphi{a}}$ and~\Ref{eq_gF} is
proven.

Equality~\Ref{eq_gF} together with Lemma~\ref{lemma_GactsonT}
completes the proof.

{\bf b)} Let $(T,(Y_1,\ldots,Y_s))$ be a path $\QS$-tableau with substitution.
Observe that for any $g\in G(\n,\g,\i)$ and a path $\al=\al_1\cdots
\al_r$ in $\Q$ we have $g\cdot X_{\al_r}\cdots
X_{\al_1}=g_{\al'}^{-1}\cdot X_{\al_r}\cdots X_{\al_1}\cdot
g_{\al''}$. Use this remark and the proof of part~a) to obtain the
claim.
\end{proof}

\begin{lemma}\label{lemma_they are invariants}
The elements from Theorem~\ref{theo_main} are invariants.
\end{lemma}
\begin{proof}
We claim that $g\cdot\Phi^{\D}(X_\be)=\Phi^{\D}(g\cdot X_\be)$ for
all $g\in G(\n,\g,\i)$ and $\be\in \Q_1^{\D}$. For $v\in \Q_0$ set  %
$$I_v=\left\{
\begin{array}{cl}
J(\n_v),&\text{if } \g_v=\Sp \\
E(\n_v),& \text{otherwise}\\
\end{array}
\right., \quad %
\de_v=\left\{
\begin{array}{rl}
-1,& \text{if }\g_v=\Sp \\
1,& \text{otherwise}\\
\end{array}
\right..
$$
Then for every $v\in \Q_0$ we have
\begin{eq}\label{eq_g_iv}
\g_{\i(v)}=\de_v I_v (g_v^{-1})^t I_v.
\end{eq}
There are two cases.
\begin{enumerate}
\item[$\bullet$] If $\be\in \Q_1$, then $g\cdot
\Phi^{\D}(X_{\be})=g\cdot X_{\be}=\Phi^{\D}(g\cdot X_{\be})$.

\item[$\bullet$] If $\be=\al^t$ for $\al\in \Q_1$, then $g\cdot
\Phi^{\D}(X_{\be})=g\cdot (I_{\al''}X_{\al}^t
I_{\al'})=I_{\al''}(g_{\al'}^{-1}X_{\al} g_{\al''})^t
I_{\al'}$. On the other hand, %
$\Phi^{\D}(g\cdot X_{\be})=\Phi^{\D}(g_{\be'}^{-1}X_{\be}
g_{\be''})=g_{\i(\al'')}^{-1}I_{\al''}X_{\al}^t I_{\al'}
g_{\i(\al')}$. Formula~\Ref{eq_g_iv} shows that $g\cdot
\Phi^{\D}(X_{\be})=\Phi^{\D}(g\cdot X_{\be})$.
\end{enumerate}
Therefore $g\cdot\Phi^{\D}(X_\be)=\Phi^{\D}(g\cdot X_\be)$ and, consequently,
$g\cdot \Phi^{\D}(f)=\Phi^{\D}(g\cdot f)$ for every $f\in
K[H(\Q^{\D},\n,\h^{\D})]$. Obviously, for a closed path
$\be_1\cdots \be_r$ in $\Q^{\D}$ and $1\leq k\leq \n_{\be_1''}$,
the element $\si_k(X_{\be_r}\cdots X_{\be_1})\in
K[H(\Q^{\D},\n,\h^{\D})]$ is a $G(\n,\g,\i)$-invariant. This
remark together with Lemma~\ref{lemma_action} concludes the proof.
\end{proof}

\section{$GL$- and $SL$-invariants of quivers}\label{section_GL_SL}
Let $\QS=(\Q,\n,\g,\h,\i)$ be a mixed quiver setting
satisfying~\Ref{eq_condition} and $G=G(\n,\g,\i)$,
$H=H(\Q,\n,\h)$. Throughout this section we assume that $\g_v$ is
$GL$ or $SL$ for all $v\in \Q_0$ and $\h_{\al}=M$ for all $\al\in
\Q_1$.

\bigskip
\noindent\textbf{ Definition (of $\QS^{\L}$).} %
Define a mixed quiver setting
$\QS^{\L}=(\Q^{\L},\n,\g,\h^{\L},\i)$ together with a mapping
$\Phi^{\L}:K[H(\Q^{\L},\n,\h^{\L})]\to K[H]$ as follows. (Here the
letter $\L$ stands for the word {\it loop}). Let $\Q_0^{\L}=\Q_0$,
$\Q_1^{\L}=\Q_1\coprod \{\al_v\,|\,v\in \Q_0,\,v<\i(v)\}$, where
$\al_v$ is a loop in the vertex $v$. For an $\al\in \Q_1$ define
$\h^{\L}_{\al}=\h_{\al}$ and $\Phi^{\L}(X_{\al})=X_{\al}$; for
$v\in \Q_0$ with $v<\i(v)$ define $\h^{\L}_{\al_v}=M$ and
$\Phi^{\L}(X_{\al_v})=E(\n_v)$.
\bigskip

It is not difficult to see that $g\cdot
\Phi^{\L}(X_{\be})=\Phi^{\L}(g\cdot X_{\be})$ for all $g\in G$ and
$\be\in \Q^{\L}_1$. Therefore
\begin{eq}\label{eq_gF_L}
g\cdot \Phi^{\L}(f)=\Phi^{\L}(g\cdot f)\text{ for all }f\in K[H(\Q^{\L},\n,\h^{\L})]
\end{eq}

\begin{theo}\label{theo_GL_SL}
Let $\QS=(\Q,\n,\g,\h,\i)$ be a mixed quiver setting
satisfying~\Ref{eq_condition}, where $\g_v$ is $GL$ or $SL$ for
all $v\in \Q_0$ and $\h_{\al}=M$ for all $\al\in \Q_1$. Then the
algebra of invariants $K[H]^{G}$ is spanned over $K$ by the
elements $\Phi^{\L}(\F_T(Y_1,\ldots,Y_s))$, where
$(T,(Y_1,\ldots,Y_s))$ is a $\QS^{\L}$-tableau with substitution of a weight
$\un{w}$ and $w_{\i(v)}=w_v$ for all $v\in \Q_0$ with $\g_v=GL$.
\end{theo}

The proof is organized as follows. The semi-invariants, i.e., the
invariants for the case of $\g_v=SL$ for all $v\in \Q_0$, has been
calculated in~\cite{LZ1}, where the case of an arbitrary quiver
setting was reduced to a {\it zigzag} quiver setting, and for
zigzag quiver settings generating systems for semi-invariants were
described. We start by rewriting these results in the language of tableaux with substitutions (see
Sections~\ref{subsection_zigzag},~\ref{subsection_SL_arbitrary}).
Afterwards it is just an exercise to reduce $GL$- and
$SL$-invariants to semi-invariants (see
Section~\ref{subsection_GL_SL}).

\subsection{Semi-invariants of zigzag quiver settings}\label{subsection_zigzag} %
The following two definitions are taken from~\cite{LZ1}.

\bigskip
\noindent\textbf{Definition.} A quiver $\Q$ is called {\it
bipartite}, if every vertex is a source (i.e. there is no arrow
ending at this vertex), or a sink (i.e. there is no arrow starting
at this vertex). A quiver setting $(\Q,\n,\g,\h,\i)$  is called a
{\it zigzag} quiver setting, if
\begin{enumerate}
\item[$\bullet$] $\Q$ is a bipartite quiver, $\g_v$ is $SL$ for
all $v\in \Q_0$ and $\h_{\al}=M$ for all $\al\in \Q_1$;

\item[$\bullet$] for every vertex $v\in \Q_0$ we have $\i(v)\neq
v$; moreover, if $v$ is a source, then $\i(v)$  is a sink and vice
versa;

\item[$\bullet$] there is no arrow $\al\in \Q_1$ with
$\al'>\i(\al')$ and $\al''>\i(\al'')$.
\end{enumerate}
\bigskip

Assume that $(\Q,\n,\g,\h,\i)$ is a zigzag quiver setting. Hence
$\Q$ can be schematically depicted as follows.
$$\begin{array}{rcccl}
l_1+l_2+1&\bullet  &           &\bullet&1 \\
&\vdots                   & \stackrel{\be(1),\ldots,\be(d_2)}{\longrightarrow}        &\vdots& \\
2l_1+ l_2&\bullet&                         &\bullet&l_1\\
                         &&  \stackrel{\al(1),\ldots,\al(d_1)}{\nearrow}   && \\
l_1+1&\bullet         &                         &\bullet&2l_1+l_2+1 \\
&\vdots                   & \stackrel{\ga(1),\ldots,\ga(d_3)}{\longrightarrow}         &\vdots& \\
l_1+l_2&\bullet       &                          &\bullet&2l_1+2l_2\\
\end{array}
$$ %
Here
\begin{itemize}
\item $\Q_0=\{1,\ldots,2l_1+2l_2\}$ and \\
$\Q_{1}=\{\al(1),\ldots,\al(d_1),\be(1),\ldots,\be(d_2),\ga(1),\ldots,\ga(d_3)\}$;

\item  the arrows $\al(1),\ldots,\al(d_1)$ go from the vertices
$l_1+1,\ldots,l_1+l_2$ to the vertices $1,\ldots,l_1$; the arrows
$\be(1),\ldots,\be(d_2)$ go from $l_1+l_2+1,\ldots,2l_1+l_2$ to
$1,\ldots,l_1$; and the arrows $\ga(1),\ldots,\ga(d_3)$ go from
$l_1+1,\ldots,l_1+l_2$ to $2l_1+l_2+1,\ldots,2l_1+2l_2$;

\item the involution $\i$ permutes vertices horizontally; consequently
$\n_{\i(v)}=\n_v$ for all $v\in\{1,\ldots,2l_1+2l_2\}$.

\end{itemize}

Fix $\un{t}=(t_1,\ldots,t_{d_1})$, $\un{r}=(r_1,\ldots,r_{d_2})$,
$\un{s}=(s_1,\ldots,s_{d_3})$ and denote by
$K[H](\un{t},\un{r},\un{s})\subset K[H]$ the space of polynomials
that have a total degree $t_k$ in variables from $X_{\al(k)}$
($1\leq k\leq d_1$), a total degree $r_k$ in variables from
$X_{\be(k)}$ ($1\leq k\leq d_2$), and a total degree $s_k$ in
variables from $X_{\ga(k)}$ ($1\leq k\leq d_3$). Further, let
$T,R,S$, respectively, be the distributions determined by $\un{t}$, $\un{r}$,
$\un{s}$, respectively, and denote $t=t_1+\cdots+t_{d_1}$,
$r=r_1+\cdots+r_{d_2}$, $s=s_1+\cdots+s_{d_3}$.

\bigskip
\noindent\textbf{Definition.}  Let $A=(A_1,\ldots,A_p)$ be a
distribution of the set $[1,t+2r]$ and let $B=(B_1,\ldots,B_q)$ be
a distribution of the set $[1,t+2s]$. The quintuple
$(\un{t},\un{r},\un{s},A,B)$ is called {\it admissible} if there
are $\un{p}=(p_1,\ldots,p_{l_1})\in \NN^{l_1}$,
$\un{q}=(q_1,\ldots,q_{l_2})\in \NN^{l_2}$ such that for
$p=p_1+\cdots+p_{l_1}$, $q=q_1+\cdots+q_{l_2}$ and for the
distributions $P,Q$, determined by $\un{p}$, $\un{q}$, we have
$\#A_j=\n_{P|j|}$ ($1\leq j\leq p$), $\#B_j=\n_{Q|j|+l_1}$ ($1\leq
j\leq q$),
$$
\bigcup_{1\leq k\leq d_1,\,  \al(k)'=i}T_k \bigcup_{1\leq k\leq
d_2,\, \be(k)'=i}(t+R_k) \bigcup_{1\leq k\leq
d_2,\,\be(k)''=\i(i)} (t+r+R_k)= \bigcup_{1\leq j\leq p,\,
P|j|=i}A_j,$$
where $1\leq i\leq l_1$, and %
$$
\bigcup_{1\leq k\leq d_1,\, \al(k)''=i}T_k \bigcup_{1\leq k\leq
d_3,\, \ga(k)'=\i(i)}(t+S_k) \bigcup_{1\leq k\leq
d_3,\,\ga(k)''=i}(t+s+S_k)=
\bigcup_{1\leq j\leq q,\, Q|j|=i-l_1}B_j,$$ %
where $l_1+1\leq i\leq l_1+l_2$. In particular,
$\sum_{i=1}^{l_1}\n_ip_i=t+2r$,
$\sum_{i=l_1+1}^{l_1+l_2}\n_iq_{i-l_1}=t+2s$. A pair
$(\un{p},\un{q})$ is called a {\it weight}.
\bigskip

The definition of the polynomial $\DP_{\un{t},\un{r},\un{s}}^{A,B}\in
K[H](\un{t},\un{r},\un{s})$ for an admissible quintuple $(\un{t},\un{r},\un{s},A,B)$
can be found in Section~5.2 of~\cite{LZ1}. This polynomial is the block partial
linearization of the $\DP$, where $\DP$ was introduced
in~\cite{LZ1} as the mixture of the determinant and the pfaffian.
(Note that $\DP$ is also a b.p.l.p., see part~4 of Example~2
from~\cite{Lop_bplp}.) By Proposition~3 from~\cite{LZ1}, in the case
$K=\QQ$ we have  %
\begin{eq}\label{eq_DP}
\DP_{\un{t},\un{r},\un{s}}^{A,B}=\frac{1}{c}\sum_{\tau_1\in
\Symmgr_A}\sum_{\tau_2\in \Symmgr_B}
F^{A,B}(\tau_1,\tau_2), %
\end{eq}
where $F^{A,B}(\tau_1,\tau_2)$ is %
$$\sign(\tau_1)\sign(\tau_2) \prod_{i=1}^t x^{T|i|}_{A\LA
\tau_1(i)\RA,   B\LA \tau_2(i)\RA} \prod_{j=1}^r y^{R|j|}_{A\LA
\tau_1(t+j)\RA, A\LA \tau_1(t+r+j)\RA} \prod_{k=1}^s
z^{S|k|}_{B\LA \tau_2(t+k)\RA, B\LA \tau_2(t+s+k)\RA}
$$
and the constant $c\in\ZZ$ depends on the quintuple. Moreover, all
coefficients in $\DP_{\un{t},\un{r},\un{s}}^{A,B}$ belong to
$\ZZ$, if the characteristic of $K$ is zero, and they belong to
$\ZZ/(\Char{K})\ZZ$, if the
characteristic of $K$ is positive. Hence using~\Ref{eq_DP}, we can define
$\DP_{\un{t},\un{r},\un{s}}^{A,B}$ over an arbitrary field. Now we
can formulate Theorem~2 of~\cite{LZ1}.
\begin{theo}\label{theo2_LZ}
If a space $K[H](\un{t},\un{r},\un{s})^{G}$ is non-zero, then the
triplet $(\un{t},\un{r},\un{s})$ is admissible. In this case
$K[H](\un{t},\un{r},\un{s})^{G}$ is spanned over $K$ by the set of
all $\DP_{\un{t},\un{r},\un{s}}^{A,B}$ for all admissible quintuples
$(\un{t},\un{r},\un{s},A,B)$.
\end{theo}

\smallskip
\noindent\textbf{Construction.} Let $(\un{t},\un{r},\un{s},A,B)$
be an admissible quintuple of a weight $(\un{p},\un{q})$. Define
$\un{n}=
(\underbrace{\n_1,\ldots,\n_1}_{p_1},\ldots,\underbrace{\n_{l_1},\ldots,\n_{l_1}}_{p_{l_1}}, %
\underbrace{\n_{l_1+1},\ldots,\n_{l_1+1}}_{q_1},\ldots, %
\underbrace{\n_{l_1+l_2},\ldots,\n_{l_1+l_2}}_{q_{l_2}})$. %
Construct a tableau with substitution $(D,(Z_1,\ldots,Z_h))$ of dimension
$\un{n}$ as follows. Arrows of $D$ are
$\{a_i,b_j,c_k\,|\,1\leq i\leq t,\, 1\leq j\leq r,\, %
1\leq k\leq s\}$, where
\begin{itemize}
\item $''a_i=A\LA i\RA$, $a_i''=A|i|$, $'a_i=B\LA i\RA$,
$a_i'=B|i|+p$;

\item $''b_j=A\LA t+j\RA$, $b_j''=A|t+j|$, $'b_j=A\LA t+r+j\RA$,
$b_j'=A|t+r+j|$;

\item $''c_k=B\LA t+k\RA$, $c_k''=B|t+k|+p$, $'c_k=B\LA t+s+k\RA$,
$c_k'=B|t+s+k|+p$.
\end{itemize}
Define $\ovphi{\cdot}$ in such a way that $\{\ovphi{a}\,|\,a\in
D\}=[1,h]$ for some $h>0$ and for any $1\leq i_1<i_2\leq t$,
$1\leq j_1<j_2\leq r$, $1\leq k_1<k_2\leq s$ we have
\begin{itemize}
\item
$\ovphi{a_{i_1}}\leq\ovphi{a_{i_2}}<\ovphi{b_{j_1}}\leq\ovphi{b_{j_2}}<
\ovphi{c_{k_1}}\leq\ovphi{c_{k_2}}$;

\item $\ovphi{a_{i_1}}=\ovphi{a_{i_2}}$ if and only if
$a_{i_1}'=a_{i_2}'$, $a_{i_1}''=a_{i_2}''$, $T|i_1|=T|i_2|$;

\item $\ovphi{b_{j_1}}=\ovphi{b_{j_2}}$ if and only if
$b_{j_1}'=b_{j_2}'$, $b_{j_1}''=b_{j_2}''$, $R|j_1|=R|j_2|$;

\item $\ovphi{c_{k_1}}=\ovphi{c_{k_2}}$ if and only if
$c_{k_1}'=c_{k_2}'$, $c_{k_1}''=c_{k_2}''$, $S|k_1|=S|k_2|$.
\end{itemize}
Define matrices $Z_1,\ldots,Z_h$ by
\begin{itemize}
\item $Z_{\ovphi{a_i}}=X_{\al(T|i|)}$,
$Z_{\ovphi{b_j}}=X_{\be(R|j|)}$, $Z_{\ovphi{c_k}}=X_{\ga(S|k|)}$
\end{itemize}
for $1\leq i\leq t,\, 1\leq j\leq r,\, %
1\leq k\leq s$. We say that $(D,(Z_1,\ldots,Z_h))$ is the tableau with substitution that
corresponds to $(\un{t},\un{r},\un{s},A,B)$.

\begin{lemma}\label{lemma_T_well}
Let $(D,(Z_1,\ldots,Z_h))$ be the tableau with substitution, of dimension $\un{n}$, that
corresponds to an admissible quintuple $(\un{t},\un{r},\un{s},A,B)$ of a weight
$(\un{p},\un{q})$. Then
\begin{enumerate}
\item[a)] $(D,(Z_1,\ldots,Z_h))$ is well defined;

\item[b)] $(D,(Z_1,\ldots,Z_h))$ is a $\QS$-tableau with substitution with the weight
$\un{w}=(p_1,\ldots,p_{l_1}$, $0,\ldots,0$, $q_1,\ldots,q_{l_2})\in\NN^{2l_1+2l_2}$;

\item[c)] $\F_D(Z_1,\ldots,Z_h)=\pm
\DP_{\un{t},\un{r},\un{s}}^{A,B}$.
\end{enumerate}
\end{lemma}
\begin{proof} {\bf a)} We claim that $X_{\al(T|i|)}$ is an
$n_{a_i''}\times n_{a_i'}$ matrix for $1\leq i\leq t$. In other
words, we should prove $n_{a_i''}=\n_{\al(T|i|)'}$ and
$n_{a_i'}=\n_{\al(T|i|)''}$. Admissibility of
$(\un{t},\un{r},\un{s},A,B)$ implies that $P|A|i||=\al(T|i|)'$ and
$Q|B|i||+l_1=\al(T|i|)''$. Thus
$$n_{a_i''}=n_{A|i|}=\n_{P|A|i||}=\n_{\al(T|i|)'}\text{ and }$$
$$n_{a_i'}=n_{B|i|+p}=\n_{Q|B|i||+l_1}=\n_{\al(T|i|)''}.$$

Similarly, we can show that $X_{\be(R|j|)}$ is an $n_{b_j''}\times
n_{b_j'}$ matrix for any $1\leq j\leq r$ and $X_{\ga(S|k|)}$ is an
$n_{c_k''}\times n_{c_k'}$ matrix for any $1\leq k\leq s$.

It is not difficult to see that each cell of $D$ is the head or
the tail of one and only one arrow. The claim follows.

{\bf b)} This part of the lemma is straightforward.

{\bf c)} Construct multi-partitions $\un{\ga}_{max}\vdash\un{t}$,
$\un{\de}_{max}\vdash\un{r}$, $\un{\la}_{max}\vdash\un{s}$, respectively, using
$(\un{t},\un{r},\un{s},A,B)$ (see part~(iii) of Definition~3
from~\cite{LZ1}) and let they determine the distributions
$\Ga_{max}$,
$\De_{max}$, $\La_{max}$, respectively. It is not difficult to see that %
$c_D=\# \Symmgr_{{\Ga}_{max}}\, \#\Symmgr_{{\De}_{max}}\, \#
\Symmgr_{{\La}_{max}}=\pm c$.  Formula~\Ref{eq_DP} concludes the
proof over the field $\QQ$ which extends to the case of an arbitrary field.
\end{proof}

The following result is a consequence of Theorem~\ref{theo2_LZ}
and Lemmas~\ref{lemma_action},~\ref{lemma_T_well}.
\begin{theo}\label{theo_zigzag}
Let $(\Q,\n,\g,\h,\i)$ be a zigzag quiver setting. Then the algebra of invariants
$K[H]^{G}$ is spanned over $K$ by the elements $\F_D(Z_1,\ldots,Z_h)$, where
$(D,(Z_1,\ldots,Z_h))$ is a $\QS$-tableau with substitution.
\end{theo}
Let us remark that under the given restrictions on $\QS$ the
generating set of Theorem~\ref{theo_zigzag} is smaller than that
of Theorem~\ref{theo_GL_SL}.

\subsection{Semi-invariants of arbitrary quivers}\label{subsection_SL_arbitrary} %
Consider a mixed quiver setting $\QS=(\Q,\n,\g,\h,\i)$ such that
$\Q$ is an arbitrary quiver, $\Q_0=\{1,\ldots,l\}$, $\g_v=SL$ for
all $v\in \Q_0$, and $\h_{\al}=M$ for all $\al\in Q_1$. Recall the construction of
the quiver $\Q^{(2)}$ from Section~4 of~\cite{LZ1}.  By
definition, $\Q_0^{(2)}=\{1,\ldots,2l\}\subset\NN$ and
$\Q_1^{(2)}=\{\ov{\al},\be_v\,|\,\al\in \Q_1,\,v\in \Q_0,\text{
and } v<\i(v)\}$, where $\be_v'=2v-1$, $\be_v''=2v$, and %
\begin{enumerate} %
\item[$\bullet$] if $\i(\al')\geq\al'$ or $\i(\al'')\geq\al''$, then $\ov{\al}'=2\al'-1$ and $\ov{\al}''=2\al''$;

\item[$\bullet$] if $\i(\al')<\al'$ and $\i(\al'')<\al''$, then $\ov{\al}'=2\i(\al'')-1$ and $\ov{\al}''=2\i(\al')$.
\end{enumerate} %
Consider the mixed quiver setting
$\QS^{(2)}=(\Q^{(2)},\n^{(2)},\g^{(2)},\h^{(2)},\i^{(2)})$, where
for all $v\in \Q_0$, $u\in \Q_0^{(2)}$, and $\ga\in \Q_0^{(2)}$ we
have %
$$\begin{array}{c}
\i^{(2)}(2v-1)=2\i(v),\quad \i^{(2)}(2v)=2\i(v)-1, \\
\n^{(2)}_{2v}=\n^{(2)}_{2v-1}=\n_v,\\
\g^{(2)}(u)=SL,\quad \h^{(2)}(\ga)=M.\\
\end{array}$$
Denote $G^{(2)}=G(\n^{(2)},\g^{(2)},\i^{(2)})$ and
$H^{(2)}=H(\Q^{(2)},\n^{(2)},\h^{(2)})$. Then $\QS^{(2)}$ is
a zigzag quiver setting and a generating system of
$K[H^{(2)}]^{G^{(2)}}$ is known. The following is the statement of Theorem~1
from~\cite{LZ1}.
\begin{theo}\label{theo1_LZ}
The homomorphism of $K$-algebras $\Phi:K[H^{(2)}]^{G^{(2)}}\to
K[H]^G$, given by
$$
\begin{array}{ccl}
\Phi(X_{\be_v})&=&E(\n_v) \text{ for } v\in \Q_0,\, v<\i(v), \\
\Phi(X_{\ov{\al}})&=&\left\{
\begin{array}{cl}
X_{\al}^t,&\text{if } \i(\al')<\al' \text{ and } \i(\al'')<\al''\\
X_{\al} ,& \text{ otherwise} \\
\end{array}
\right. \\%
\end{array}
$$
for $\al\in \Q_1$, is a surjective mapping.
\end{theo}

\begin{lemma}\label{lemma_tableauL}
For every $\QS^{(2)}$-tableau with substitution $(D,(Z_1,\ldots,Z_{s}))$ there is a
$\QS^{\L}$-tableau with substitution $(T,(Y_1,\ldots,Y_s))$ such that
$\Phi^{\L}(\F_T(Y_1,\ldots,Y_s))=\Phi(\F_D(Z_1,\ldots,Z_{s}))$.
\end{lemma}
\begin{proof}
Let $(D,(Z_1,\ldots,Z_{s}))$ be $\QS^{(2)}$-tableau with substitution of weight
$\un{w}=(w_1,\ldots,w_{2l})$ and dimension $\un{n}$ and let $W$ be the distribution
determined by $\un{w}$.

For every $v\in \Q_0$ the vertex $2v$ of $\Q^{(2)}$ is a source and $\i^{(2)}(2v)$ is a
sink. Assume that $w_{2v}\neq0$. Hence there is an $a\in D$ such that $a'$ or $a''$ lies
in $W_{2v}$. The definition of a $\QS^{(2)}$-tableau with substitution implies a
contradiction. Thus $w_{2v}=0$ for all $v\in \Q_0$.

Let $\un{u}=(u_1,\ldots,u_l)$, where $u_v=w_{2v-1}$ for all
$v\in \Q_0$, and let $U$ be the distribution determined by
$\un{u}$. Note that for any $1\leq i\leq w_1+\cdots+w_{2l}$ we have
$W|i|=2U|i|-1$.

Define a tableau with substitution $(T,(Y_1,\ldots,Y_s))$ of dimension $\un{n}$ as
follows. Define $T=\{b_a\,|\,a\in D\}$ and for $a\in D$ define $\ovphi{b_a}=\ovphi{a}$.

Consider $\ga\in \Q_1^{(2)}$
such that $Z_{\ovphi{a}}=X_{\ga}$ and $W|a'|=\i^{(2)}(\ga'')$,
$W|a''|=\ga'$.

If $\ga=\ov{\al}$ for an $\al\in \Q_1$, then define
$Y_{\ovphi{a}}=X_{\al}$. If $\i(\al')<\al'$ and $\i(\al'')<\al''$,
then $U|a'|=\al'$, $U|a''|=\i(\al'')$ and we define ${b_a}''=a'$,
${b_a}'=a''$, $''{b_a}={}'a$, $'{b_a}={}''a$. If $\i(\al')>\al'$
or $\i(\al'')>\al''$, then $U|a'|=\i(\al'')$, $U|a''|=\al'$ and we define $b_a=a$.

If $\ga=\be_v$ for $v\in \Q_0$, $v<\i(v)$, then $U|a'|=\i(\al_v'')$, $U|a''|=\al_v'$
for $\al_v\in \Q_0^{\L}$ and we define $Y_{\ovphi{a}}=X_{\al_v}$, $b_a=a$.

This describes a $\QS^{\L}$-tableau with substitution $(T,(Y_1,\ldots,Y_s))$ of dimension
$\un{n}$ and weight $\un{u}$. Obviously,
$\Phi^{\L}(\F_T(Y_1,\ldots,Y_s))=\Phi(\F_D(Z_1,\ldots,Z_{s}))$.
\end{proof}

\begin{theo}\label{theo_SL}
Let $\QS=(\Q,\n,\g,\h,\i)$ be a mixed quiver setting satisfying~\Ref{eq_condition},
$\g_v=SL$ for all $v\in \Q_0$, and $\h_{\al}=M$ for all $\al\in \Q_1$. Then the algebra
of invariants $K[H]^{G}$ is spanned over $K$ by the elements
$\Phi^{\L}(\F_T(Y_1,\ldots,Y_s))$, where $(T,(Y_1,\ldots,Y_s))$ is a $\QS^{\L}$-tableau
with substitution.
\end{theo}
\begin{proof}
Theorems~\ref{theo_zigzag},~\ref{theo1_LZ} together with
Lemma~\ref{lemma_tableauL} show that the invariants belong to
$K$-span of the elements from the theorem. Formula~\Ref{eq_gF_L}
and Lemma~\ref{lemma_action} completes the proof.
\end{proof}

\subsection{Proof of Theorem~2}\label{subsection_GL_SL}
Let $\Q_0=\{1,\ldots,l\}$. Define
$\g^{(1)}=(\g^{(1)}_1,\ldots,\g^{(1)}_l)$ by $\g^{(1)}_v=SL$ for
all $v\in \Q_0$ and denote $G^{(1)}=G(\n,\g^{(1)},\i)$. Since
$G^{(1)}\subset G$, we have $K[H]^{G}\subset K[H]^{G^{(1)}}$. Also
note that $\Q_0^{\L}=\Q_0$.

Let $f\in K[H]^G$ be a polynomial of a multidegree $\un{t}=(t_{\al})_{\al\in \Q_1}$. By
Theorem~\ref{theo_SL}, $f=\sum_j\la_j \Phi^{\L}(\F_{T_j}(Y_{j,1},\ldots,Y_{j,s_j}))$,
where $\la_j\in K$, $(T_j,(Y_{j,1},\ldots,Y_{j,s_j}))$ is a $\QS^{\L}$-tableau with
substitution of a weight $\un{w}_j=(w_{j,1},\ldots,w_{j,l})$, and the multidegree of
$\Phi^{\L}(\F_{T_j}(Y_{j,1},\ldots,Y_{j,s_j}))$ is $\un{t}$. Lemma~\ref{lemma_action}
together with formula~\Ref{eq_gF_L} and restriction~\Ref{eq_condition} imply that for any
$g\in G$ we have
$$g\cdot f=\sum_j\la_j \Phi^{\L}(\F_{T_j}(Y_{j,1},\ldots,Y_{j,s_j})) \left(
\prod_{v\in \Q_0,\, v<\i(v)}
\det(g_v)^{w_{j,\i(v)}-w_{j,v}}\right).$$ %
It is not difficult to see that $w_{j,\i(v)}-w_{j,v}$ does not
depend on $j$, therefore we denote it by $\mu_v$. Hence,
$$g\cdot f=\prod_{v\in \Q_0,\, v<\i(v)}\det(g_v)^{\mu_v}\cdot
f=f,$$ %
and the statement of the theorem follows immediately.

\section{Reduction}\label{section_reduction}
Consider an arbitrary mixed quiver setting $\QS=(\Q,\n,\g,\h,\i)$
satisfying~\Ref{eq_condition}. As usual, define
$G=G(\n,\g,\i)$ and $H=H(\Q,\n,\h)$.

\bigskip
\noindent\textbf{Definition (of $\QS^{\R}$).} %
Let a mixed quiver setting
$\QS^{\R}=(\Q^{\R},\n^{\R},\g^{\R},\h^{\R},\i^{\R})$ and a mapping
$\Phi^{\R}:K[H(\Q^{\R},\n^{\R},\h^{\R})]\to K[H]$ be a result of
the following procedure. (Here the letter $R$ stands for the word
{\it reduction}). Denote indeterminates of
$K[H(\Q^{\R},\n^{\R},\h^{\R})]$ by $\ov{x}_{ij}^\be$, where
$\be\in\Q_1^{\R}$, and denote the corresponding generic matrix by
$\ov{X}_{\be}=(\ov{x}_{ij}^\be)$. Initially assume that $\QS^{\R}$
is $\QS$ and for every $\al\in \Q_1$ we have
$\Phi^{\R}(\ov{X}_\al)=X_\al$. Afterwards we change $\QS^{\R}$
and $\Phi^{\R}$ by performing steps~a)--c).
\begin{enumerate}
\item[a)] For every $v\in \Q_0$ such that $\g_v$ is $O$ or $\Sp$
we add a new vertex $\ov{v}$ to $\Q^{\R}$ and set
$\i^{\R}(v)=\ov{v}>v$, $\n_{\ov{v}}^{\R}=\n_v$,
$\g^{\R}_{v}=\g^{\R}_{\ov{v}}=GL$. Also we add two new arrows $\be_v$,
$\ga_v$ to $\Q^{\R}$ such that $\be_v'=v$, $\be_v''=\ov{v}$,
$\ga_v'=\ov{v}$, $\ga_v''=v$.  By definition, if $\g_{v}=O$, then
$\Phi^{\R}(\ov{X}_{\be_v})=\Phi^{\R}(\ov{X}_{\ga_v})=E(\n_v)$,
otherwise $\Phi^{\R}(\ov{X}_{\be_v})=
\Phi^{\R}(\ov{X}_{\ga_v})=J(\n_v)$.

\item[b)] For every $v\in \Q_0$ with $\g_v=SO$ we add a new vertex
$\ov{v}$ to $\Q^{\R}$ and set $\i^{\R}(v)=\ov{v}>v$,
$\n_{\ov{v}}^{\R}=\n_v$, $\g^{\R}_v=\g^{\R}_{\ov{v}}=SL$. Also we add a
new arrow $\be_v$ to $\Q^{\R}$ such that $\be_v'=v$,
$\be_v''=\ov{v}$.  By definition,
$\Phi^{\R}(\ov{X}_{\be_v})=E(\n_v)$;

\item[c)] Define $\h^{\R}_{\be}=M$ for all $\be\in\Q_1^{R}$.
\end{enumerate}
Note that $\Q_0\subset \Q_0^{\R}$, $\Q_1\subset \Q_1^{\R}$. For
short, we write $G^{\R}$ for $G(\n^{\R},\g^{\R},\i^{\R})$ and
$H^{\R}$ for $H(\Q^{\R},\n^{\R},\h^{\R})$.
\bigskip

Let $\psi:G\to G^{\R}$ be the natural embedding, where
$\psi(g)_v=g_v$ for $v\in \Q_0$, and
$\psi(g)_{\ov{v}}=(g_v^{-1})^t$ for $v\in \Q_0$ with
$\g_v\in\{O,\Sp,SO\}$. It is not difficult to see that $g\cdot
\Phi^{\R}(\ov{X}_{\be})=\Phi^{\R}(\psi(g)\cdot \ov{X}_{\be})$ for
all $g\in G$ and $\be\in \Q^{\R}_1$. Hence
\begin{eq}\label{eq_gF_R}
g\cdot \Phi^{\R}(f)=\Phi^{\R}(\psi(g)\cdot f)
\end{eq}
for all $f\in K[H^{\R}]$. Thus $\Phi^{\R}(f)$ lies in $K[H]^G$ for
all $f$ from $K[H^{\R}]^{G^{\R}}$.

\begin{theo}\label{theo_surjection}
The restriction $\Phi^{\R}:K[H^{\R}]^{G^{\R}}\to K[H]^G$ is a
surjective mapping.
\end{theo}

\noindent{}Since a generating system for $K[H^{\R}]^{G^{\R}}$ is
known (Theorem~\ref{theo_GL_SL}), Theorem~\ref{theo_surjection}
gives us a generating system for $K[H]^G$. In order to prove
this theorem we need some facts from the theory of modules with
good filtrations.

\subsection{Good filtrations}\label{subsection_filtration}
Consider an affine algebraic group $G$ with the coordinate ring
$K[G]$, a closed subgroup $B$ of $G$, and a rational $B$-module
$A$. The tensor product $K[G]\otimes A$ over $K$ is naturally a
rational $G\times B$ module with respect to the action $(g,b)\cdot
(f\otimes a)=f^{(g,b)}\otimes ba$ for $g\in G$, $b\in B$, $f\in
K[G]$, $a\in A$, and $f^{(g,b)}(x)=f(g^{-1}xb)$. The set of
$B$-fixed points of $K[G]\otimes A$ is a rational $G$-module,
called the {\it induced} module $\Ind_{B}^G A=(K[G]\otimes A)^B$
(see, for example,~\cite{Grosshans97}).

By a {\it good filtration} of a rational $G$-module $M$ we mean an
ascending chain $0=M_0\subseteq M_1\subseteq \cdots \subseteq M$
of submodules with $\cup_{i\in\NN} M_i=M$ and for $i>0$ the module
$M_{i+1}/M_i$ is either zero or is induced from a one-dimensional
module over a fixed Borel subgroup of $G$.
\begin{enumerate}
\item[$\bullet$] An affine $G$-variety $H$ is called {\it good} if
its coordinate ring $K[H]$ has a good filtration, where, as usual,
$G$ acts on $K[H]$ by the rule $(g\cdot f)(x)=f(g^{-1}\cdot x)$ for $g\in G$,
$f\in K[H]$, and $h\in H$.

\item[$\bullet$] A pair of affine $G$-varieties $(H_2,H_1)$ such that $H_1$ is
a closed subvariety of $H_2$ given by an ideal $I$, i.e.,
$K[H_1]=K[H_2]/I$, is called a {\it good $G$-pair}, if $H_2$
is good, and $I$ has a good filtration. %
\end{enumerate}

We list some standard properties of modules with good filtrations
(see~\cite{Donkin85},~\cite{Mathieu90}).
\begin{theo}\begin{enumerate}\label{theo_good_modules}
\item[a)] If $$0\to L\to M\to N\to 0$$ is a short exact sequence
of $G$-modules and $L$ has a good filtration, then %
$$0\to L^G\to M^G\to N^G\to 0$$
is exact.

\item[b)] If $M\subseteq N$ are modules with good filtrations,
then $N/M$ also has a good filtration.

\item[c)] If $M$, $N$ are $G$-modules with good filtrations, then
$M\otimes N$, considered under the diagonal action
of $G$, also has a good filtration.
\end{enumerate}
\end{theo}

\begin{cor}\label{cor_good_pairs}
If $(H_2,H_1)$ is a good $G$-pair, then $H_1$ is a good $G$-module
and the mapping $K[H_2]^G\to K[H_1]^G$, induced by the natural
surjection $K[H_2]\to K[H_1]$, is also a surjection.
\end{cor}

The following two lemmas will help us to construct good pairs.
\begin{lemma}\label{lemma_good_pairs_general}
\begin{enumerate}
\item[a)] Let $M$ be a $GL(n)$-module. Then $M$ has a good
$GL(n)$-filtration if and only if $M$ has a good
$SL(n)$-filtration.

\item[b)] Let $0=M_0\subseteq M_1\subseteq \cdots \subseteq M_r=M$
be an ascending chain of submodules such that $M_{i+1}/M_i$ has a
good filtration for any $1\leq i<r$. Then $M$ has a good
filtration.

\item[c)] Let $M$ be a $G_1\times G_2$-module such that $G_2$ acts
trivially on $M$. Then $M$ has a good $G_1\times G_2$-filtration
if and only if $M$ has a good $G_1$-filtration.

\item[d)] If $(H_2,N)$ and $(N,H_1)$ are good pairs, then
$(H_2,H_1)$ is a good pair.

\item[e)] If $(H_1,N_1)$ and $(H_2,N_2)$ are good $G$-pairs, then
$(H_1\oplus H_2,N_1\oplus N_2)$ is a good $G$-pair, provided that $G$ acts
on the direct sums diagonally.
\end{enumerate}
\end{lemma}
\begin{proof} For parts~b),~c),~d), respectively, see Proposition~1.2a~(iv),
Proposition~1.2e~(ii), Lemma~1.3a~(i) of~\cite{Donkin90},
respectively. Since $SL(n)$ is the derived subgroup of $GL(n)$,
part~a) follows from part~(i) of Lemma~1.4 of~\cite{Donkin93}.

Now we prove part~e). Part~c) of Theorem~\ref{theo_good_modules} implies
that $K[H_1\oplus H_2]=K[H_1]\otimes K[H_2]$ has a good
filtration. We have $K[N_i]=K[H_i]/I_i$ for
$i=1,2$ and some ideals $I_i$. Consider the short exact sequence %
$$0\to I \to K[H_1]\otimes K[H_2]\to
K[N_1]\otimes K[N_2]\to 0,$$ %
where $I=I_1\otimes K[H_2] + K[H_1]\otimes I_2$. %
Since $I\,/\,{I_1\otimes I_2}=(I_1\otimes K[N_2])\, \oplus\,
(K[N_1]\otimes I_2)$ has a good filtration and $I_1\otimes I_2$
has a good filtration (see Theorem~\ref{theo_good_modules},~c)),
then $I$ has a good filtration by part~b).
\end{proof}

Denote by $S^{+}_0(n)$ and $S^{+}_1(n)$, respectively, subsets of $S^{+}(n)$ consisting
of invertible matrices and of matrices with the
determinant $1$, respectively. Similarly denote by $S^{-}_0(n)$
the subset of $S^{-}(n)$ that consists of invertible matrices. Let
$\{1, \ast\}$ be the group of order two, where the symbol $1$
stands for the identity element of the group. Given $\be\in\{1,
\ast\}$, a vector space $V=K^n$, and $g\in GL(n)$, we write
$$V^{\be}=\left\{
\begin{array}{cl}
V,& \be=1\\
V^{\ast},&\be=\ast\\
\end{array}
\right.,\quad %
g^{\be}=\left\{
\begin{array}{ccc}
g,&\be=1\\
(g^{-1})^t,& \be=\ast\\
\end{array}
\right..$$ %

\begin{lemma}\label{lemma_good_pairs}
\begin{enumerate}%
\item[a)] Let $GL(n)$ act on $K^{n\times n}$ by the formula $g\cdot A=gAg^t$ for
$g\in GL(n)$, $A\in K^{n\times n}$. Then
$(K^{n\times n},S^{+}(n))$, %
$(K^{n\times n},S^{-}(n))$, %
$(GL(n),S^{+}_0(n))$, and %
$(GL(n),S^{-}_0(n))$ %
are good $GL(n)$-pairs.

\item[b)] Let the group $SL(n)$ act on $K^{n\times n}$ in the same
manner as in part~a). Then $(S^{+}(n),S^{+}_1(n))$ is a good
$SL(n)$-pair.

\item[c)] Let $GL(n)$ act on $H=K^{n\times n}\oplus K^{n\times n}$
by the formula $g\cdot (A_1,A_2)=(gA_1g^t,(g^{-1})^tA_2 g^{-1})$ for $g\in
GL(n)$ and $(A_1,A_2)\in H$. Then $(H,\{(A_1,A_2)\in
H\,|\,A_1A_2=E\})$ is a good $GL(n)$-pair.

\item[d)] Let $V=(K^n)^{\be}$ and $W=(K^m)^{\ga}$ for
$\be,\ga\in\{1,\ast\}$. If $GL(n)\times GL(m)$ acts on
$M=K^{m\times n}\simeq \Hom_K(V,W)$ in the natural way, i.e.,
$(g_1,g_2)\cdot A=g_2^{\ga}A(g_1^{\be})^{-1}$, then $M$ is a good
$GL(n)\times GL(m)$-module.
\end{enumerate}
\end{lemma}
\begin{proof}
For part~a) see Lemma~1.3 of~\cite{Zubkov99} and its proof. For
part~c) see the reasoning at the end of
Section~1 of~\cite{Zubkov99}.

To prove part~b) let $K[S^{+}(n)]=K[x_{ij}\,|\,1\leq i,j\leq
n]$, where $x_{ij}=x_{ji}$ and consider an $n\times n$ matrix
$X=(x_{ij})$ and an ideal $I\vartriangleleft K[S^{+}(n)]$
generated by $\det(X)-1$.
Obviously $K[S^{+}_1(n)]=K[S^{+}(n)]/I$. %
A mapping $f\mapsto f\cdot (\det(X)-1)$ is an isomorphism of
$SL(n)$-modules $K[S^{+}(n)]$ and $I$. This fact together with
part~a) of this lemma, part~a) of
Lemma~\ref{lemma_good_pairs_general} and
Corollary~\ref{cor_good_pairs} imply that $S^{+}(n)$ is a good
$SL(n)$-module and $I$ has a good $SL(n)$-filtration.

To prove part~d) notice that $K[M]\simeq S(W^{\ast}\otimes V)$ as
$GL(n)\times GL(m)$-modules (for example, see Lemma~1 of~\cite{LZ1}).
Hence it is enough to show that $S^r(W^{\ast}\otimes V)$ has a
good filtration for every $r>0$. Consider ABW-filtration for
$S^r(W^{\ast}\otimes V)$ (see Section~1.3 of~\cite{ZubkovI} or
Theorem~2.4 of~\cite{DZ01}). Its factors are
$L_{\la}(W^{\ast})\otimes L_{\la}(V)$, where $L_{\la}(V)$ denotes a
Schur module. The module $L_{\la}(V)$ has a good
$GL(n)$-filtration and $L_{\la}(W^{\ast})$ has a good
$GL(m)$-filtration. By part~c) of
Lemma~\ref{lemma_good_pairs_general} both of these modules have
good $GL(n)\times GL(m)$-filtrations. Part~c) of
Theorem~\ref{theo_good_modules} concludes the proof.
\end{proof}

We have the following version of {\it Frobenius reciprocity}.
\begin{lemma}\label{lemma_Frobenius}
(c.f. Lemma~8.1 of~\cite{Grosshans97}) If $G$ is a closed subgroup of an algebraic
group $C$ and $H$ is an affine $C$-variety (i.e. $C$ acts rationally on $H$),
then the algebra of invariants
$(K[H]\otimes K[C/G])^{C}$ is isomorphic to $K[H]^G$, where $C$
acts on $C/G$ by left multiplication and it acts on the tensor
product diagonally. The isomorphism is given by the mapping $a\otimes f\mapsto
f(1_{C/G})a$ for $a\in K[H]$ and $f\in K[C/G]$.
\end{lemma}

\noindent{}Frobenius reciprocity enables us to reduce
investigation of the invariants of $G$ to invariants of a bigger algebraic group
$C$ such that it is determined by smaller number of equations than $G$ and it is more simple
to calculate its invariants than invariants of $G$.

\subsection{Proof of Theorem~7}

We split the proof into several lemmas. Let us recall that
$H=H(\Q,\n,\h)=\bigoplus_{\al\in\Q_1}H_{\al} $, where
\begin{enumerate}
\item[a)] $H_{\al}=K^{\n_{\al'}\times\n_{\al''}}$, if $\h_{\al}=M$;

\item[b)] $H_{\al}=S^{+}(\n_{\al'})$ and $\g_{\al'}$ is $O$ or $SO$, if $\al'=\al''$ and $\h_{\al}=S^{+}$;

\item[c)] $H_{\al}=L^{+}(\n_{\al'})$ and $\g_{\al'}=\Sp$, if $\al'=\al''$ and $\h_{\al}=L^{+}$;

\item[d)] $H_{\al}=S^{+}(\n_{\al'})$, $\i(\al')=\al''$, and $\g_{\al'}$ is
$GL$ or $SL$, if $\al'\neq\al''$ and $\h_{\al}=S^{+}$;

\item[e)] analogues statements are true if we substitute $S^{-}$ for $S^{+}$ in parts~b),~d) and if we substitute
$L^{-}$ for $L^{+}$ in part~c).
\end{enumerate}

Denote %
$$H_1= %
\bigoplus_{\al\in\Q_1,\,\h_{\al}\neq L^{+},L^{-}} H_{\al} %
\bigoplus_{\al\in\Q_1,\,\h_{\al}=L^{+}}S^{+}(\n_{\al'})
\bigoplus_{\al\in\Q_1,\,\h_{\al}=L^{-}}S^{-}(\n_{\al'}).
$$ %
and define the action of $G$ on $H_1$ by %
\begin{eq}\label{eq_action_G_on_H_1}
(g\cdot h)_{\al}=\left\{
\begin{array}{cl}
g_{\al'}h_{\al}g_{\al'}^{t},& \text{ if } \h_{\al}\neq M\text{ and }\al\text{ is a loop}\\
g_{\al'}h_{\al}g_{\al''}^{-1},& \text{ otherwise }\\
\end{array}
\right.,
\end{eq}%
where $g=(g_v)_{v\in\Q_0}$ lies in $G$, $h=(h_{\al})_{\al\in\Q_1}$ lies in
$H_1$, and $\al\in\Q_1$. As a consequence of the following remark the spaces $H$
and $H_1$ are isomorphic as $G$-modules.

\begin{remark}\label{remark_1null}
Let $g\in \Sp(n)$, $A\in L^{+}(n)$, and $B\in S^{+}(n)$. Define
the action of $\Sp(n)$ on $L^{+}(n)$ by $g\cdot A=gAg^{-1}$ and
its action on $S^{+}(n)$ by $g\cdot B=gBg^t$. Then the mapping
$L^{+}(n)\to S^{+}(n)$ given by $A\mapsto AJ$ is an isomorphism of
$\Sp(n)$-modules. The assertion remains valid if we replace $L^{+}$ and $S^{+}$, respectively,
by $L^{-}$ and $S^{-}$, respectively.
\end{remark}

Rewrite $G$ as $G=P\times \prod G_v$, where $v$ ranges over
vertices of $\Q$ such that $\g_v\in\{O,\Sp,SO\}$. The group $G$ is
a subgroup of
$$C=P\times
\prod_{v\in\Q_0,\,\g_v=O,\Sp}GL(n_v)
\prod_{v\in\Q_0,\,\g_v=SO}SL(n_v).$$ %
Notice that $H_1$ is also a $C$-module, where the action is defined
by~\Ref{eq_action_G_on_H_1}. Frobenius reciprocity (see
Lemma~\ref{lemma_Frobenius}) gives the isomorphism
$$(K[H_1]\otimes K[C/G])^C\simeq K[H_1]^G.$$

Consider the $C$-module %
$$W_1=
\bigoplus_{v\in\Q_0,\,\g_v=O}S^{+}_0(\n_v)
\bigoplus_{v\in\Q_0,\,\g_v=\Sp}S^{-}_0(\n_v)
\bigoplus_{v\in\Q_0,\,\g_v=SO}S^{+}_1(\n_v),
$$%
where $S_0^{+}(n)$, $S_0^{-}(n)$, and $S_1^{+}(n)$ were defined
just before Lemma~\ref{lemma_good_pairs}. Here the action is given
by
\begin{eq}\label{eq_action_C_on_W_1}
(g\cdot w)_v=g_v w_v g_v^t\text{ for all }v\in\Q_0\text{ with
}\,\g_v\in\{O,\Sp,SO\},
\end{eq}%
where $g\in C$ and $w=(w_v)_{v\in\Q_0,\,\g_v\in\{O,\Sp,SO\}}$
belongs to $W_1$. The mapping $\pi_1:C/G\to W_1$ given by
$$\pi_1(\ov{g})=\left\{
\begin{array}{cl}
g_vg_v^t,& \text{if }g_v\text{ is }O\text{ or }SO\\
g_vJg_v^t,& \text{if }g_v=\Sp\\
\end{array}
\right.
$$
($\ov{g}$ is the equivalence class of $g\in C$, $v\in \Q_0$, and
$\g_v\in\{O,\Sp,SO\}$) is well defined.

\begin{lemma}\label{lemma_W_1}
The mapping $\pi_1$ is a $C$-equivariant isomorphism of the
algebraic groups $C/G$ and $W_1$.
\end{lemma}
\begin{proof} Mappings
$$GL(n)/O(n)\to S^{+}_0(n),\,\,\ov{g}\mapsto gg^t \text{ and}$$
$$GL(n)/\Sp(n)\to S^{-}_0(n),\,\,\ov{g}\mapsto gJg^t$$
are $GL(n)$-equivariant isomorphisms of the algebraic groups,
where $GL(n)$ acts on the left hand sides by left multiplication
and it acts on the right hand side by the formula $g\cdot A=gAg^t$ (see
Lemma~1.2 of~\cite{Zubkov99}). Similarly, a mapping
$$SL(n)/SO(n)\to S^{+}_1(n),\,\,\ov{g}\mapsto gg^t,$$
is an $SL(n)$-equivariant isomorphisms of the algebraic groups,
where the actions are the same as above. This completes the proof.
\end{proof}

The space $H_1$ is contained in
$$H_2=\bigoplus_{\al\in\Q_1} K^{\n_{\al'}\times \n_{\al''}}$$
and the space $W_1$ is contained in
$$W_2=
\bigoplus_{v\in\Q_0,\,\g_v\text{ is }O\text{ or }\Sp}GL(\n_v)
\bigoplus_{v\in\Q_0,\,\g_v=SO}S^{+}(\n_v).$$ %
Define the action of $C$ on $H_2$ and $W_2$, respectively, by the
formulas~\Ref{eq_action_G_on_H_1} and~\Ref{eq_action_C_on_W_1},
respectively.

Lemmas~\ref{lemma_good_pairs_general},~\ref{lemma_good_pairs}
imply the following result.

\begin{lemma}\label{lemma_good_1}
The $C$-pair $(H_2\oplus W_2,H_1\oplus W_1)$ is a good one.
\end{lemma}

We refer to the elements of the space
$$W_3=
\bigoplus_{v\in\Q_0,\,\g_v=O,\Sp} %
(K^{\n_v\times\n_v}\oplus K^{\n_v\times\n_v})    %
\bigoplus_{v\in\Q_0,\,\g_v=SO}K^{\n_v\times \n_v}$$ %
as
$w=(w_v,w_{\ov{u}})_{v,u\in\Q_0,\,\g_v\in\{O,\Sp,SO\},\,\g_u\in\{O,\Sp\}}$.
Endow $W_3$ with the structure of $C$-module by
$$(g\cdot w)_v=g_vw_vg_v^t\text{ for all }v\in\Q_0\text{ with }\g_v\in\{O,\Sp,SO\},$$
$$(g\cdot w)_{\ov{u}}=(g_u^{-1})^t w_{\ov{u}}g_u^{-1}\text{ for all }u\in\Q_0\text{ with }\g_u\in\{O,\Sp\}$$
for $g\in C$. Define the embedding
$$\pi_2:W_2\to W_3$$
of $C$-spaces by $\pi_2(w)_v=w_v$ and $\pi_2(w)_{\ov{u}}=w_u^{-1}$.

The following lemma is a consequence of
Lemmas~\ref{lemma_good_pairs_general},~\ref{lemma_good_pairs}.

\begin{lemma}\label{lemma_good_2}
The $C$-pair $(H_2\oplus W_3,H_2\oplus \pi_2(W_2))$ is a good one.
\end{lemma}

Using Corollary~\ref{cor_good_pairs} together with
Lemmas~\ref{lemma_W_1},~\ref{lemma_good_1},~\ref{lemma_good_2},
and Frobenius reciprocity (see Lemma~\ref{lemma_Frobenius}) we
obtain a surjection
$$\Phi:K[H_2\oplus W_3]^C\to K[H]^G.$$ %

Consider the mixed quiver setting $\QS^{\R}$. Remove each loop
$\al\in\Q_1^{\R}$ with $\h_{\al}\neq M$ and replace it by a new arrow $\be_{\al}$ with  $\be_{\al}'=v$,
$\be_{\al}''=\ov{v}$, where $\al'=\al''=v$. Denote the
resulting quiver by $\Q^{\New}$ and set
$\QS^{\New}=(\Q^{\New},\n^{\R},\g^{\R},\h^{\New},\i^{\R})$, where
$\h^{\New}_{\be}=M$ for all $\be\in \Q_1^{\New}$. For short, we
write $G^{\New}$ for $G(\n^{\R},\g^{\R},\i^{\R})$ and write
$H^{\New}$ for $H(\Q^{\New},\n^{\R},\h^{\New})$. Denote
indeterminates of $K[H^{\New}]$ by $\ov{x}_{ij}^\be$, where
$\be\in\Q_1^{\New}$, and denote the corresponding generic matrix
by $\ov{X}_{\be}=(\ov{x}_{ij}^\be)$. Then $C=G^{\New}$ and
$G^{\New}$-modules $H_2\oplus W_3$ and $H^{\New}$ are equal.

Define a mapping
$\Phi^{\New}:K[H^{\New}]\to K[H]$ by %
$$\Phi^{\New}(\ov{X}_{\ga})=\left\{
\begin{array}{cl}
X_{\al},&\text{if }\ga=\be_{\al}\text{ for a loop }\al\in\Q_1
\text{ with }\h_{\al}\in\{S^{+},S^{-}\}\\
X_{\al}J(\n_v),&\text{if }\ga=\be_{\al}
\text{ for a loop }\al\in\Q_1\text{ with }\h_{\al}\in\{L^{+},L^{-}\}\\
\Phi^{\R}(\ov{X}_{\be}),&\text{if }\ga\in \Q_1^{\R}\cap
\Q_1^{\New}
\end{array}
\right. $$ %
for $\ga\in\Q_1^{\New}$. It is not difficult to see that the
restriction of $\Phi^{\New}$ to $G^{\New}$-invariants coincides
with $\Phi$.

Let $\F_T(Y_1,\ldots,Y_s)$, where $(T,(Y_1,\ldots,Y_s))$ is a $\QS^{\New}$-tableau with
substitution of a weight $\un{w}$, be a $G^{\New}$-invariant. Then there is a path
$\QS^{\R}$-tableau with substitution
$(D,(Z_1,\ldots,Z_s))$ such that %
$$\Phi^{\New}(\F_T(Y_1,\ldots,Y_s))=\Phi^{\R}(\F_D(Z_1,\ldots,Z_s))$$
and the weight of $(D,(Z_1,\ldots,Z_s))$ is $\un{w}$. (Note that
$\Q_0^{\New}=\Q_0^R$.) Lemma~\ref{lemma_action} implies that
$\F_D(Z_1,\ldots,Z_s)$ is a $G^{\R}$-invariant.
Theorem~\ref{theo_GL_SL} together with~\Ref{eq_gF_R} shows that
$$\Phi^{\New}(K[H^{\New}]^{G^{\New}})\subset
\Phi^{\R}(K[H^{\R}]^{G^{\R}})\subset K[H]^G.$$ %
The fact that $\Phi$ is surjective completes the proof of
Theorem~\ref{theo_surjection}.

\section{Path $\QS$-tableaux with substitutions and good $\QS$-tableaux with substitutions}\label{section_tableau_pairs}

Consider a tableau with substitution $(T,(X_1,\ldots,X_s))$ of dimension $\un{n}\in\NN^m$
and numbers $1\leq q_1<q_2\leq m$ with $n_{q_1}=n_{q_2}$. Let us recall some concepts
introduced in~\cite{Lop_bplp}.

Denote by $\M$ the monoid freely generated by letters
$x_1,x_2,\ldots$, $x_1^t,x_2^t,\ldots$  Let $\M_T$ be the
submonoid of $\M$, generated by $x_1,\ldots,x_s$,
$x_1^t,\ldots,x_s^t$. For short, we will write
$1,\ldots,s,1^t,\ldots,s^t$ instead of
$x_1,\ldots,x_s,x_1^t,\ldots,x_s^t$. Given $a\in T$, we consider
$\ovphi{a}\in\{1,\ldots,s\}$ as an element of $\M_T$.

Given $u\in\M_T$, define the matrix $X_u$ by the following rules:
\begin{enumerate}
\item[$\bullet$] $X_{j^t}=X_j^t$ for any $1\leq j\leq s$;

\item[$\bullet$] $X_{vw}=\left\{%
\begin{array}{cl}
X_v X_w,&\text{if the product of these matrices is well defined} \\
0,& \text{otherwise} \\
\end{array}
\right.$ %
\\for $v,w\in \M_T$.
\end{enumerate}

For an arrow $a\in T$ denote by $a^t$ the {\it transpose arrow},
i.e., by definition $(a^t)''=a'$, $(a^t)'=a''$, $''(a^t)={}'a$,
$'(a^t)={}''a$, $\ovphi{a^t}=\ovphi{a}^t\in \M_T$. Obviously,
$(a^t)^t=a$.

We write $a\stackrel{t}{\in} T$ if $a\in T$ or $a^t\in T$.

\bigskip
\noindent{\bf Definition (of paths in $T$).} We say that
$a_1,a_2\stackrel{t}{\in}T$ are {\it successive} in $T$ (with respect to $q_1$ and $q_2$),  if
$a_1',a_2''\in\{q_1,q_2\}$, $a_1'\neq a_2''$, $'a_1={}''a_2$.

A word $a=a_1\cdots a_r$ with
$a_1,\ldots,a_r\stackrel{t}{\in}T$, is called a {\it path} in $T$
with respect to columns $q_1$ and $q_2$, if $a_i$, $a_{i+1}$ are
successive (with respect to $q_1$ and $q_2$) for any $1\leq i\leq r-1$. In this case by definition
$\ovphi{a}=\ovphi{a_1} \cdots \ovphi{a_r}\in \M_T$ and
$a^t=a_r^t\cdots a_1^t$ is a path in $T$; we denote
$a_r',{}'a_r,a_1'',{}''a_1$, respectively, by $a',{}'a,a'',{}''a$, respectively.
If the columns $q_1$ and $q_2$ are fixed, we refer to $a$ as a path in $T$.
A path $a_1\cdots a_r$ is {\it closed} if $a_r,a_1$ are
successive;   in particular,  $a_1'',a_r'\in\{q_1,q_2\}$.

Given $\xi\in \Symmgr_{n_{q_1}}$ we will also use the notations
$(T^{\xi},(Y_1,\ldots,Y_s))$ and $(\widetilde{T}^{\xi},Y_{\widetilde{T}^{\xi}})$,
respectively, for the tableaux with substitutions of dimensions $\un{n}$ and $\un{d}$,
respectively, defined in Section~4 of~\cite{Lop_bplp}. Here
$Y_{\widetilde{T}^{\xi}}=(Z_1,\ldots,Z_h)$ for some matrices $Z_1,\ldots,Z_h$.

\bigskip
Throughout this section we assume that $\QS=(\Q,\n,\g,\h,\i)$ is a
mixed quiver setting satisfying~\Ref{eq_condition} and
\begin{eq}\label{eq_not_Sp}
\g_v\neq \Sp \text{ for all }v\in \Q_0.
\end{eq}%
Consider an arrow $\be\in\Q_1^{\D}$. If $\be\in\Q_1$ and $\h_{\be}=M$, then
$\be^t\in\Q_1^{\D}$ is given by the definition of a mixed double quiver setting.
Otherwise, we set
$$\be^t=\left\{
\begin{array}{cl}
\al,& \text{if }\be=\al^t,\,\h_{\al}=M, \text{ and }\al\in\Q_1\\
\be,& \text{if }\be\in\Q_1 \text{ and }\h_{\be}\neq M\\
\end{array}
\right..$$%
If $\be=\be_1\cdots\be_r$ is a path in $\Q^{\D}$, denote by
$\be^t=\be_r^t\cdots\be_1^t$ a path in $\Q^{\D}$.
Condition~\Ref{eq_not_Sp} implies that
$\Phi^{\D}(X_{\be_1^t}\cdots
X_{\be_r^t})=\Phi^{\D}(X_{\be_r}\cdots X_{\be_1})^t$.

\bigskip
\noindent\textbf{Definition (of good $\QS$-tableaux with substitutions).} A tableau with
substitution $(T,(Y_1,\ldots,Y_s))$ of dimension $\un{n}\in \NN^m$ is called a {\it good
$\QS$-tableau with substitution}, if for some {\it weight}
$\un{w}=(w_1,\ldots,w_l)\in\NN^l$ and the distribution $W$, determined by $\un{w}$, one
has
\begin{enumerate}
\item[$\bullet$]
$\un{n}=(\underbrace{\n_1,\ldots,\n_1}_{w_1},\ldots,\underbrace{\n_l\ldots,\n_l}_{w_l})$;

\item[$\bullet$] if $a\in T$, then there exists a path
$\be=\be_1\cdots \be_r$ in $\Q^{\D}$ (where $\be_1,\ldots,\be_r\in
\Q_1^{\D}$) such that $Y_{\ovphi{a}}=\Phi^{\D}(X_{\be_r}\cdots
X_{\be_1})$, $W|a'|=\i(\be'')$, and $W|a''|=\be'$.
\end{enumerate}

\begin{remark}\label{remark_good_and_path_tabl}
\begin{enumerate}
\item[a) ] Any path $\QS$-tableau with substitution is also a good $\QS$-tableau with
substitution.

\item[b) ] For every good $\QS$-tableau with substitution $(T,(Y_1,\ldots,Y_s))$ of a
weight $\un{w}$ there is a path $\QS^{\D}$-tableau with substitution
$(T,(Z_1,\ldots,Z_s))$ of weight $\un{w}$ such that for every $a\in T$ we have the
equality $Z_{\ovphi{a}}=X_{\be_r}\cdots X_{\be_1}$ for a path $\be_1\cdots\be_r$ from the
definition. Moreover, $\F_{T}(Y_1,\ldots,Y_s)=\Phi^{\D}(\F_{T}(Z_1,\ldots,Z_s))$.
\end{enumerate}
\end{remark}

\begin{lemma}\label{lemma_Q_good_tableau}
Consider a good $\QS$-tableau with substitution $(T,(Y_1,\ldots,Y_s))$ of a weight
$\un{w}$, that determines the distribution $W$, and let $\un{n}\in\NN^{m}$ be the
dimension of $(T,(Y_1,\ldots,Y_s))$. Assume that numbers $1\leq q_1<q_2\leq m$ satisfy
$\i(W|q_1|)=W|q_2|$ and $b=b_1\cdots b_p$ (where $b_1,\ldots,b_p\stackrel{t}{\in}T$) is a
path in $T$ with respect to columns $q_1$ and $q_2$.

Then there is a path $\be=\be_1\cdots \be_r$ in $\Q^{\D}$ (where
$\be_1,\ldots,\be_r\in \Q_1^{\D}$) such that
$Y_{\ovphi{b}}=\Phi^{\D}(X_{\be_r}\cdots X_{\be_1})$,
$W|b'|=\i(\be'')$, and $W|b''|=\be'$. In particular, if $b$ is a
closed path in $T$, then $\be$ is a closed path in $\Q^{\D}$.
\end{lemma}
\begin{proof}
Obviously, it is enough to prove the lemma for $p=2$. For $j\in\{1,2\}$, the condition
$b_j\stackrel{t}{\in} T$ implies that there is an $a_j\in T$ with $b_j\in\{a_j,a_j^t\}$.
By the definition of a good $\QS$-tableau with substitution, there is a path
$\al_j=\al_{j,1}\cdots \al_{j,r_j}$ in $\Q^{\D}$ (where $\al_{j,1},\ldots,\al_{j,r_j}\in
\Q_1^{\D}$) such that $Y_{\ovphi{a_j}}=\Phi^{\D}(X_{\al_{j,r_j}}\cdots X_{\al_{j,1}})$,
$W|a'_j|=\i(\al''_j)$, and $W|a''_j|=\al'_j$.

Assume $b_1=a_1$ and $b_2=a_2^t$. Then $\{q_1,q_2\}=\{a_1',a_2'\}$.
Define
$\be=\al_2^t\al_1=\al_{2,r_2}^t\cdots\al_{2,1}^t\al_{1,1}\cdots\al_{1,r_1}=\be_1\cdots\be_r$,
where $\be_1,\ldots,\be_r\in \Q_1^{\D}$, $r=r_1+r_2$. Since
$(\al_2^t)'=\i(\al_2'')=W|a_2'|=\i(W|a_1'|)=\al_1''$, we infer that $\be$ is a
path in $\Q^{\D}$. Obviously,
$$Y_{\ovphi{b}}=Y_{\ovphi{a_1}} Y_{\ovphi{a_2}}^t=%
\Phi^{\D}(X_{\al_{1,r_1}}\cdots X_{\al_{1,1}} X_{\al_{2,1}}^t
\cdots X_{\al_{2,r_2}}^t)=\Phi^{\D}(X_{\be_r}\cdots X_{\be_1}).
$$
By the definitions of $\al_1$, $\al_2$, we have %
$W|b'|=W|b_2'|=W|a_2''|=\al_2'=\i(\be'')$
and %
$W|b''|=W|b_1''|=W|a_1''|=\al_1'=\be'$. This proves the claim
for the given case.

The remaining cases can be treated analogously.
\end{proof}

\begin{theo}\label{theo_decomposition}
Let $\QS=(\Q,\n,\g,\h,\i)$ be a mixed quiver setting satisfying
conditions~\Ref{eq_condition} and~\Ref{eq_not_Sp}. Let $(T,(Y_1,\ldots,Y_s))$ be a path
$\QS$-tableau with substitution of a weight $\un{w}$, that determines the distribution
$W$. Then $\F_T(Y_1,\ldots,Y_s)$ is a polynomial over $K$ in
$\Phi^{\D}(\si_k(X_{\be_r}\cdots X_{\be_1}))$ and $\Phi^{\D}(\F_D(Z_1,\ldots,Z_h))$,
where
\begin{enumerate}
\item[1.] $\be_1\cdots \be_r$ ranges over closed paths in $\Q^{\D}$ and
$1\leq k\leq \n_{\be_1''}$;       

\item[2.] $(D,(Z_1,\ldots,Z_h))$ ranges over path $\QS^{\D}$-tableaux with substitutions
of a weight $\un{u}$ such that            

\begin{enumerate}
\item[a)] for every $v\in \Q_0$ we have $u_v\leq w_v$ and
$w_{\i(v)}-w_v=u_{\i(v)}-u_v$;   

\item[b)] if $\g_v$ is $GL$ or $SL$ for $v\in \Q_0$, then
$u_{\i(v)}=0$ or $u_{v}=0$;      

\item[c)] if $\g_v$ is $O$ or $SO$ for $v\in \Q_0$, then
$u_v\leq 1$ and $\i(v)=v$.       
\end{enumerate}
\end{enumerate}
\end{theo}
\begin{proof} We have that $(T,(Y_1,\ldots,Y_s))$ is a good $\QS$-tableau with substitution
of weight $\un{w}$. If conditions~b) and~c) are valid for $\un{u}=\un{w}$, then the
statement is trivial. Otherwise, there is a $v\in \Q_0$ such that $\i(v)\neq i$,
$w_{\i(v)}\neq 0$, $w_v\neq0$ or $\i(v)=v$, $w_v\geq2$. Consider $1\leq q_1<q_2\leq m$
such that $W|q_1|=v$, $W|q_2|=\i(v)$. In particular, $\i(W|q_1|)=W|q_2|$ and
$n_{q_1}=\n_{v}=\n_{\i(v)}=n_{q_2}$. Apply the decomposition formula (see Theorem~2
from~\cite{Lop_bplp}) to the $q_1$-th and the $q_2$-th columns of $T$. Then
$\F_T(Y_1,\ldots,Y_s)$ is a
polynomial in $\si_k(Y_{\ovphi{c}})$, $\F_D(Z_1,\ldots,Z_h)$, where %
\begin{itemize}
\item $c$ is a closed path in $T^{\xi}$ for some $\xi\in
\Symmgr_{n_{q_1}}$;

\item $(D,(Z_1,\ldots,Z_h))\!=\!(\widetilde{T}^{\pi},Y_{\widetilde{T}^{\pi}})$ is a
tableau with substitution of dimension $\un{d}$ for some $\pi\in \Symmgr_{n_{q_1}}$.
\end{itemize}
Here $\un{d}\in\NN^{m-2}$ is obtained from $\un{n}$ by eliminating the $q_1$-th and
$q_2$-th coordinates. Applying Lemma~\ref{lemma_Q_good_tableau} to paths of $T^{\pi}$, we
obtain that $(\widetilde{T}^{\pi},Y_{\widetilde{T}^{\pi}})$ is a good $\QS$-tableau with
substitution. Lemma~\ref{lemma_Q_good_tableau} also imply that
$\si_k(Y_{\ovphi{c}})=\si_k(\Phi^{\D}(X_{\be_r}\cdots X_{\be_1}))=
\Phi^{\D}(\si_k(X_{\be_r}\cdots X_{\be_1})$ for a closed path $\be_1\cdots \be_r$ in
$\Q^{\D}$.

Repeat this procedure for $(D,(Z_1,\ldots,Z_h))$ and so on. Finally we see that
$\F_T(Y_1,\ldots,Y_s)$ is a polynomial in $\Phi^{\D}(\si_k(X_{\be_r}\cdots X_{\be_1}))$
and $\F_D(Z_1,\ldots,Z_h)$ such that condition~1 is valid and $(D,(Z_1,\ldots,Z_h))$ is a
good $\QS$-tableau with substitution of a weight $\un{u}$ satisfying~a),~b), and~c).
Part~b) of Remark~\ref{remark_good_and_path_tabl} completes the proof.
\end{proof}

\section{Proof of Theorem~1}\label{section_proof_main}
Consider a mixed quiver setting $\QS=(\Q,\n,\g,\h,\i)$
satisfying~\Ref{eq_condition}. As usual, define $G=G(\n,\g,\i)$,
$H=H(\Q,\n,\h)$. In Lemma~\ref{lemma_they
are invariants} we have shown that the elements from
Theorem~\ref{theo_main} are invariants.

Successive applications of Theorem~\ref{theo_surjection} and Theorem~\ref{theo_GL_SL}
yield the generating system that consists of images of $\F_T(Y_1,\ldots,Y_s)$ for some
tableaux with substitutions $(T,(Y_1,\ldots,Y_s))$. Application of
Theorem~\ref{theo_decomposition} to $\F_T(Y_1,\ldots,Y_s)$ gives the mixed quiver
setting, which we denote by $\QS^0=(\Q^0,\n^0,\g^0,\h^0,\i^0)$, and the mapping
$\Phi^0:K[H^0]\to K[H]$, where $H^0$ stands for $H(\Q^0,\n^0,\h^0)$, such that
$$\Q_0^0=\Q_0^{\R}=\Q_0\coprod \{\ov{v}\,|\,v\in\Q_0,\,\g_v\in\{O,\Sp,SO\}\},$$%
$$\Q_1^0=\Q_1 \coprod \{\al^t\,|\,\al\in Q_0\}%
\coprod \{\al_v,\al_v^t\,|\,v\in\Q_0,v\leq\i(v)\}$$
$$\coprod \{\be_v,\be_v^t\,|\,v\in\Q_0,\g_v\in\{O,\Sp,SO\}\}
\coprod \{\ga_v,\ga_v^t\,|\,v\in\Q_0,\g_v\in\{O,\Sp\}\},$$ %
$$\i^0(v)=\left\{
\begin{array}{cl}
\i(v),&\text{ if }\g_v\in\{GL,SL\} \\
\ov{v},&\text{ if }\g_v\in\{O,\Sp,SO\}    \\
\end{array}
\right.,\quad\text{ where }v\in\Q_0,
$$
$$\g^0_v=\g^0_{\i^0(v)}=\left\{
\begin{array}{cl}
GL,&\text{if }\g_v\in\{GL,O,\Sp\} \\
SL,&\text{if }\g_v\in\{SL,SO\}    \\
\end{array}
\right.,\quad\text{ where }v\in\Q_0,$$ %
$$\h_{\be}^0=M \text{ for all }\be\in\Q_1^0, \text{ and
}\n^0=\n^{\R}.$$%
Here $\be_v,\be_v^t$ go from $\i^0(v)=\ov{v}$ to $v$, arrows
$\ga_v,\ga_v^t$ go in the opposite direction (i.e. from $v$ to $\ov{v}$), and $\al_v, \al_v^t$
are loops in $v,\i^0(v)$, respectively. Denote indeterminates of
$K[H^0]$ by $\ov{x}_{ij}^\be$ for $\be\in\Q_1^0$ and denote
the corresponding generic matrix by
$\ov{X}_{\be}=(\ov{x}_{ij}^\be)$. The mapping $\Phi^0$ is
determined by
$$\Phi^0(\ov{X}_{\al_v})=E,$$ %
$$\Phi^0(\ov{X}_{\be_v})=\Phi^0(\ov{X}_{\ga_v})=\left\{
\begin{array}{cl}
E,&\text{if }\g_v\text{ is }O\text{ or }SO\\
J,&\text{if }\g_v=\Sp\\
\end{array}
\right.,$$ %
$$\Phi^0(\ov{X}_{\al})=X_{\al}\text{ for }\al\in Q_1,$$ %
$$\Phi^0(\ov{X}_{\be^t})=X_{\be}^t\text{ for }\be\in Q_1^0.$$ %
The algebra of invariants $K[H]^G$ is generated by the elements
$\Phi^0(\F_T(Y_1,\ldots,Y_s))$ and $\Phi^0(\si_k(\ov{X}_{\ga_r}\cdots
\ov{X}_{\ga_1}))$, where
\begin{itemize}
\item $(T,(Y_1,\ldots,Y_s))$ is a path $\QS^0$-tableau with substitution of a weight
$\un{u}$, satisfying the conditions
\begin{enumerate}
\item[a)] if $\g^0_v=GL$ ($v\in \Q_0^0$), then $u_{\i^0(v)}=u_v=0$,
\item[b)] if $\g^0_v=SL$ ($v\in \Q_0^0$), then $u_{\i^0(v)}=0$ or
$u_v=0$.
\end{enumerate}

\item $\ga_1\cdots \ga_r$ is a closed path in $\Q^0$, $1\leq k\leq
\n^0_{\ga_1''}$.
\end{itemize}

\begin{lemma}\label{lemma_X}
Let $\ga=\ga_1\cdots \ga_r$ be a path in $\Q^0$ and
$X=\Phi^0(\ov{X}_{\ga_r}\cdots \ov{X}_{\ga_1})$.
\begin{enumerate}
\item[a)] If $X$ is not a matrix over $K$, then there is a path
$\be=\be_1\cdots \be_p$ in $\Q^{\D}$ such that
$X=\pm\Phi^{\D}(X_{\be_p}\cdots X_{\be_1})$ and
$$\be'=\left\{
\begin{array}{cl}
\ga',& \text{if }\ga'\in \Q_0\\
\i^0(\ga'),& \text{if } \ga'\notin \Q_0\\
\end{array}
\right.\,\,,\quad %
\be''=\left\{
\begin{array}{cl}
\ga'',& \text{if } \ga''\in \Q_0\\
\i^0(\ga''),& \text{if } \ga''\notin \Q_0\\
\end{array}
\right..$$ %
In particular, if $\ga$ is closed, then $\be$ is closed.

\item[b)] If $X$ is a matrix over $K$, then there is a $v\in\Q_0$
such that $\ga_i',\ga_i''\in\{v,\i^0(v)\}$ for any $1\leq i\leq
r$. Moreover,

if $\g_v$ is $GL$ or $SL$, then for any $1\leq i\leq r$ an arrow
$\ga_i$ is a loop;

if $\g_v\neq\Sp$, then $X=\pm E$, otherwise, $X$ is $\pm E$ or
$\pm J$.
\end{enumerate}
\end{lemma}
\begin{proof}
\textbf{a)} Eliminate arrows $\al_v,\al_v^t$ from the path $\ga$
and obtain a path in $\Q^0$. Then eliminate arrows
$\be_v,\be_v^t,\ga_v,\ga_v^t$ to get the required path. To prove
it, see the following pictures, where we depicted some vertex
$v\in\Q_0$ with $\g_v\in\{O,\Sp,SO\}$ and two arrows
$\al_1,\al_2\in\Q_1$.
$$
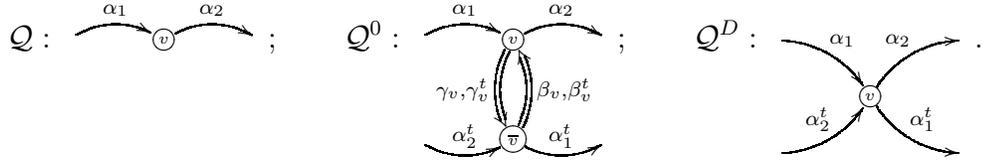
\Q:\xymatrix@C=1cm@R=1cm{ %
\ar@/^/@{->}[r]^{\al_1}&\vtx{v}\ar@/^/@{->}[r]^{\al_2}&\\
}; %
\qquad
\Q^0:\xymatrix@C=1cm@R=1cm{ %
\ar@/^/@{->}[r]^{\al_1}&\vtx{v}
\ar@2@/^/@{<-}[d]^{\be_v,\be_v^t}
\ar@2@/_/@{->}[d]_{\ga_v,\ga_v^t}
\ar@/^/@{->}[r]^{\al_2}&\\
\ar@/_/@{->}[r]^{\al_2^t}&\vtx{\ov{v}}\ar@/_/@{->}[r]^{\al_1^t}&\\
}; %
\qquad
\Q^{\D}:\xymatrix@C=1cm@R=0.5cm{ %
\ar@/^/@{->}[rd]^{\al_1}&&\\
&\vtx{v}\ar@/^/@{->}[ru]^{\al_2}\ar@/_/@{->}[rd]^{\al_1^t}&\\
\ar@/_/@{->}[ru]^{\al_2^t}&&\\
}. %
$$

\textbf{b)} This case is trivial.
\end{proof}

The last lemma shows that we can reduce $\QS^0$ to a mixed quiver
setting $\QS^1=(\Q^1,\n^{\D},\g^{\D},\h^1,\i^{\D})$, which almost
coincides with $\QS^{\D}$. To do so, define $\Q_0^1=\Q_0^{\D}=\Q_0$
and $\Q_1^1=\Q_1^{\D}\coprod \{\de_1,\ldots,\de_q\}$, where
$\de_i$ is a loop and $q\geq0$. Then consider a mapping
$\Phi^1:K[H(\Q^0,\n^0,\h^0)]\to K[H]$ such that
$$\text{if }\be\in\Q_0,\text{ then }\Phi^1(X_{\be})=\Phi^{\D}(X_{\be}),
\text{ otherwise }\Phi^1(X_{\be})\text{ is }\pm E \text{ or }\pm
J,$$%
and assume that for every path $\QS^0$-tableau with substitution $(T,(Y_1,\ldots,Y_s))$
satisfying the property formulated before Lemma~\ref{lemma_X} there is a path
$\QS^1$-tableau with substitution $(T,(Z_1,\ldots,Z_s))$ satisfying
\begin{enumerate}
\item[a)] $\Phi^0(\F_T(Y_1,\ldots,Y_s))=\pm
\Phi^1(\F_T(Z_1,\ldots,Z_s))$;

\item[b)] if the weight of $(T,(Y_1,\ldots,Y_s))$ is $\un{u}$,
then the weight of $(T,(Z_1,\ldots,Z_s))$ is $\un{w}$, where for
every $v\in\Q_0$ we have%
$$w_v=\left\{
\begin{array}{cl}
u_v,&\text{if }\g_v\text{ is }GL\text{ or }SL\\
u_v+u_{\ov{v}},&\text{if }\g_v\in\{O,\Sp,SO\}\\
\end{array}
\right..$$
\end{enumerate}
Consider a vertex $v\in\Q_0$. There are three possibilities.
\begin{enumerate}
\item[$\bullet$] If $\g_v\in\{GL,O,\Sp\}$, then equalities
$w_{\i(v)}=w_v=0$ imply that $\de_i$ is not a loop in $v$ for all
$i$.

\item[$\bullet$] If $\g_v$ is $SL$, then, since $w_v=0$ or
$w_{\i(v)}=0$, we see that $\de_i$ is not a loop in $v$ for all
$i$.

\item[$\bullet$] Let $\g_v$ be $SO$ and let $\de_i$ be a loop in
$v$ for some $i$. Then $\Phi^1(X_{\de_i})=\pm E$.
Theorem~\ref{theo_decomposition} together with the facts

\begin{enumerate}
\item[a)] $\F_T(Y_1,\ldots,Y_s)=0$ for any tableau with substitution
$(T,(Y_1,\ldots,Y_s))$ with $Y_{\ovphi{a}}=E$ and $a'=a''$ for some $a\in T$;

\item[b)] every $\QS^{\D\D}$-tableau with substitution, where $\QS^{\D\D}$ is a mixed
double quiver setting of $\QS^{\D}$, is also a path $\QS^{\D}$-tableau with substitution;
\end{enumerate}
completes the proof.
\end{enumerate}


\bigskip
\noindent{\bf Acknowledgements.} This paper was written during
author's visit to University of Antwerp, sponsored by Marie Curie
Research Training Network Liegrits. The author is grateful for this support. The author
would also like to thank Fred Van Oystaeyen for hospitality and
Alexander Zubkov for his generous help in proving of
Theorem~\ref{theo_surjection}.


\end{document}